# A Deterministic and Linear Model of Dynamic Optimization

**Somdeb Lahiri**

**ORCID: https://orcid.org/0000-0002-5247-3497**

**(Formerly with) PD Energy University, Gandhinagar (EU-G), India.**

**Abstract**

We introduce a model of infinite horizon linear dynamic optimization and obtain results concerning existence of solution and satisfaction of the "competitive condition" and "transversality condition" being unconditionally sufficient for optimality of a trajectory. We also show that under some mild restrictions the optimal trajectory satisfies the "Euler condition" and a related "transversality condition". The optimal trajectory satisfies the functional equation of dynamic programming. Under an additional convexity assumption for the two-period constraint sets, we show that the optimal value function is concave and continuous. Linearity "bites" when it comes to the definition of optimal decision rules which can no longer be guaranteed to be single-valued. We show that if all the two-period constraint sets are convex, then the optimal decision rule is an upper semi-continuous correspondence. For linear cake-eating problems, we obtain monotonicity results for the optimal value function and a "conditional monotonicity" result for optimal decision rules. We also introduce the concept of a two-phase linear cake eating problem and obtain a necessary condition that must be satisfied by all solutions of such problems. We show that for a class of linear dynamic optimization problems, known as interlinked linear dynamic optimization problems, a slightly modified version of the functional equation of dynamic programming is satisfied.

## 1. Introduction:

Notwithstanding the fact that programming in linear spaces- both finite dimensional and infinite dimensional- begins with Hurwicz (1958), in economics it has been sustained almost entirely by "infinite horizon dynamic optimization".

The problem we are concerned with here, is a generalized version of a problem which is widely known as the "cake-eating problem" (see Schmidheiny and Walti (2002)). In the cake-eating problem, a decision-maker faces the problem of choosing a piece from a given cake of finite size/weight over an infinite time horizon. At time $t = 0$, the decision-maker starts with the entire cake. At each point in time $t = 1, 0, \ldots$, the decision-maker chooses a piece from the cake that is available at that point in time and consumes it instantaneously, thereby determining the quantity of the cake that is available for future consumption. Consuming a piece of the cake gives the decision-maker instantaneous satisfaction (that may be time-dependent). Starting with a cake of a predetermined dimension, at each point of time where a decision has to be made, the decision-maker chooses the size of the cake that will remain for future consumption, thereby determining the size of the piece of cake that the decision-maker consumes at that point of time.

Let $X = [0, b] \subset \mathbb{R}$ (the set of real numbers), with $b > 0$ denote the **set of available alternatives**. In the context of the cake-eating problem, the initial size of the cake could be any number in the interval $[0, b]$. With $\mathbb{N}$ denoting the set of natural number (i.e., the set of strictly positive integers) let $\mathbb{N}^0$ denote $\mathbb{N} \cup \{0\}$, i.e., the set of non-negative integers. Time is measured in discrete periods $t\in\mathbb{N}^0$. At each time 't' an alternative is chosen, and the chosen alternative is denoted by $x_t\in X$. In the context of the cake-eating problem, $x_0$ would denote the initial size of the cake and for $t\in\mathbb{N}^0$, $x_t - x_{t+1}$ would denote the size of the piece of the cake that the decision-maker chooses to consume in period t, thereby determining the size of the cake (i.e., $x_{t+1}$) that would be available in period t+1.

A very general form of the typical dynamic optimization problem along the lines defined in section 5.1 of Sorger (2015) is the following:

Given $x\in X$: Maximize $\sum_{t=0}^{\infty} u_t(x_t, x_{t+1})$, subject to the infinite sequence $<x_t|t\in\mathbb{N}^0>$ satisfying the constraints, $(x_t, x_{t+1})\in\Omega_t$ for all $t\in\mathbb{N}^0$ and $x_0 = x\in X$, where for all $t\in\mathbb{N}^0$, $u_t$: $X\times X\rightarrow \mathbb{R}$ is the **utility function at time-period** t, $\Omega_t \subset X\times X$ is the **two-period constraint set at time-period** t, and $x\in X$ is the **initial choice**.

We will refer to the problem defined above as **the general dynamic optimization problem**.

In the context of the cake-eating problem, for all $t\in\mathbb{N}^0$, $\Omega_t = \{(x_t, x_{t+1})\in X\times X|\ x_{t+1}\in[0, x_t]\}$.

**Note 1:** The exact mathematical interpretation of the expression (formula) $\sum_{t=0}^{\infty} u_t(x_t, x_{t+1})$ is $\lim_{T\to\infty}(\sum_{t=0}^{T} u_t(x_t, x_{t+1}))$. Thus, the problem we are concerned with here is in the domain of asymptotic analysis, which is very different from infinite dimensional analysis, the latter being-in the context of social sciences-no more than an artifact.

While our discussion here partially parallels (to be precise "is based on") the first five sections of the introduction to dynamic optimization theory available in Mitra (2000), it is important to pin-point the similarities and differences between the framework of analysis in Mitra (2000) and the framework of analysis that we are concerned with here. First of all, Mitra (2000) assumes that (i) there exists a real number $\delta$ in the open interval (0, 1) and a function u: $X\times X\rightarrow \mathbb{R}$ such that for all $x, y\in X$ and $t\in\mathbb{N}^0$, $u_t(x, y) = \delta^t u(x, y)$; (ii) there exists $\Omega$ such that $\Omega = \Omega_t$ for all $t\in\mathbb{N}^0$. In our framework, the objective function is linear, with the coefficients determined by an infinite sequence of time-dependent real numbers. Thus, while the most general form of the objective function in Mitra (2000) allows for non-linearity, our objective function does not. However, in our framework, the time-dependent coefficients of the linear objective function are not restricted in any way at all. Linearity implies continuity and therefore issues such as upper semi-continuous utility functions versus continuous utility functions, that Mitra (2000) devotes considerable space to, find no mention in our analysis.

Further, Mitra (2000) assumes throughout the discussion reported there, that the "two-period constraint sets" are time invariant. We make no such assumption, except that a consequence of the cake-eating problem that we discuss in section 6 is the fact that two-period constraint sets are time invariant.

The "justification" for the linearity assumption we invoke is proposition 2.1 in this analysis, is an "asymptotic version" of a result for differentiable concave functions in Euclidean spaces that is available in Lahiri (2024) and some references therein. The result says, that for an optimization problem with a differentiable concave objective function and linear constraints, a necessary and sufficient condition for an alternative to be an optimal solution is that the alternative solves the problem with the same linear constraints and a linear objective function whose coefficient for each variable is the partial derivative of the objective function with respect to the variable evaluated at the alternative.

In the context of the model of linear dynamic optimization that we propose, one possible interpretation of the coefficients of the objective function could be associated with uncertainty about the duration of the optimization problem and instantaneous satisfaction in time-period t being proportional to $x_t - x_{t+1}$ for all t. There is also a linear version of the "optimal growth" model of Kamihigashi and Roy (2006), that can be formulated as a linear dynamic optimization problem, and we discuss this formulation in section 3. We refer to this linear model as "linear model of optimal wealth accumulation". We also introduce "interlinked" linear dynamic optimization problems, where instantaneous pay-offs at any time period depends on the decision variable for the same time period as well as the decision variable for the immediately next time-period.

For the general version of the linear dynamic optimization problem, our first result in section 4 says that an optimal solution from each and every initial point always exists. Our second result in the same section says that under certain mild restrictions, an optimal trajectory satisfies the "Euler condition" (similar to the one in section 5.2 of Sorger (2015)) and the "transversality condition". Our third result in this section says that if a trajectory from a given initial point satisfies the "competitive condition" and the trajectory satisfies a related "transversality condition", then it is an optimal solution. This result is very similar to proposition 4.1 in Mitra (2000).

Section 5 is about dynamic programming in the context of our model. The first proposition asserts that if the two-period constraint sets are convex in every time period, then the optimal value function is concave and continuous. The second proposition in this section is that a trajectory is a solution if and only if it satisfies the "functional equation of dynamic programming". It is worth emphasizing that the satisfaction of the "functional equation of dynamic programming" does not require any convexity assumption about any two-period constraint set.

Section 6 is about optimal decision rules and since linearity of the objective function invalidates the concept of "strict concavity", we need to allow for the possibility of the optimal decision rules being multivalued. Further, since in general no kind of "stationarity" is implied in our framework, we have to allow for the possibility of the optimal decision rule being time-dependent. We show that if all the two-period constraint sets are convex, then the optimal decision rule at any period, is an "upper-semicontinuous correspondence" in the sense of Debreu (1959) and satisfies certain properties that automatically follow from its definition.

Section 7 is about cake-eating problems. The interesting point about a cake-eating problem is that given an initial alternative, any trajectory starting from there is "feasible" if and only if it is "non-increasing". Our first proposition in this section says that at any time-period, if the coefficient of the objective function at that time-period is non-negative, then the optimal value function is non-increasing, and if the coefficient is positive, then the optimal value function is strictly increasing. The corresponding results in the case of non-positive and strictly negative coefficients, are conditional. The second proposition in section 7 says, that if all time-dependent optimal value functions are monotonically non-decreasing, then never consuming the cake is an optimal trajectory. The third proposition in this section says, that if at a time-period more is available, then for any alternative chosen for the next period by the optimal decision rule at the greater initial alternative if it is "no greater" than the lesser initial alternative, and any alternative chosen by the optimal decision rule at the lesser initial alternative but is unchosen for the greater initial alternative, the alternative chosen for the greater initial alternative is greater than the alternative chosen for the lesser initial alternative. This is a "conditional monotonicity" result for the optimal decision rules associated with cake-eating problems. We also introduce the concept of a two-phase linear cake-eating problem in which all coefficients of the objective function are positive up to a certain time period $T^+$ and thereafter from a period $T^-$ all coefficients of the objective function are negative. We show that for any solution of such a problem, the entire cake is available for consumption till the beginning of time-period $T^+$ and there is no cake left to be consumed beginning with time-period $T^-$. The entire cake is consumed between time periods $T^+$ and $T^-$ - 1.

In section 8 we consider a very important type of linear dynamic optimization problem, which we refer to as "interlinked linear dynamic optimization problem" and whose major characteristic we have already described in an earlier paragraph. In this case, the major theoretical result is a slightly modified version of the functional equation of dynamic programming. Once again, no convexity assumption is required for any two-period constraint set.

While Mitra (2000) and Sorger (2015) are cited several times in this paper- and justly so- it is imperative to point out a crucial difference of perspective between our work here and those of the other authors we mentioned. Mitra (2000) is obviously concerned with dynamic optimization, only to the extent that it is applicable to "optimal growth" and related issues in economics. Sorger (2015) relies on macroeconomic applications and their implications to motivate dynamic optimization techniques in chapter 5 of the book. We, on the other hand, are concerned with infinite horizon "linear" dynamic optimization. More importantly, time dependence in a non-trivial sense and linearity are the two issues that are of utmost importance to us in this paper.

**2. Framework of Analysis and Mathematical Motivation:**

Given a positive integer n and n-tuples $\chi = (x_1, \ldots, x_n)$, $\eta = (y_1, \ldots, y_n) \in \mathbb{R}^n$, the **inner product** of $\chi$ and $\eta$ denoted by $<\chi, \eta> = \sum_{j=1}^{n} \chi_j \eta_j$ and the **(Euclidean-) norm** of $\chi$ denoted by $||\chi||_2 = \sqrt[2]{<\chi, \chi>}$.

We assume the following for what follows in this section:

(M1) (i) For all $t\in\mathbb{N}^0$, $\Omega_t$ is a non-empty, closed and convex subset of X×X; (ii) For all $t\in\mathbb{N}^0$ and $x\in X$, there exists $y\in X$ (possibly depending on x) such that $(x, y)\in\Omega_t$.

(M2) For all $t\in\mathbb{N}^0$, $u_t$: X×X→ ℝ is bounded, concave and continuously differentiable, where for all (x, y)∈ X×X, $Du_t(x,y) = (D_1u_t(x,y), D_2u_t(x,y)) = (\lim_{h\to 0, x+h\in X} \frac{u_t(x+h,y)-u_t(x,y)}{h}$, $\lim_{h\to 0, y+h\in X} \frac{u_t(x,y+h)-u_t(x,y)}{h})$.

(M3) For all sequences $<x_t|t\in\mathbb{N}^0>$ and $<y_t|t\in\mathbb{N}^0>$ with $(x_t, y_t)\in X\times X$ for all $t\in\mathbb{N}^0$: $\sum_{t=0}^{\infty} |D_i u_t(x_t, y_t)| < +\infty$, $i\in\{1, 2\}$.

**Important notation:** For all $t\in\mathbb{N}^0$ and $x\in X$, let $\Omega_t(x) = \{y|\ (x, y)\in \Omega_t\}$

**Note 2:** The requirement that for all $t\in\mathbb{N}^0$, $\Omega_t$ is a non-empty, closed and convex subset of X×X and the assumption that X = [0, b] implies that for all $t\in\mathbb{N}^0$ and $x\in X$, $\Omega_t(x)$ is a non-empty closed interval in [0, b], though possibly a singleton (degenerate).

**Note 3:** Part (ii) of assumption M1 says that for all $t\in\mathbb{N}^0$, $\{x\in X|\Omega_t(x) \neq \phi\} = X$. If for $t\in\mathbb{N}$, and x, y, z, w∈X, (x, z), (y, w)∈ $\Omega_t$ then the convexity of $\Omega_t$ assumed in part (i) of M1, implies that for all θ∈(0, 1) $(\theta x + (1-\theta)y, \theta z + (1-\theta)w)\in\Omega_t$.

**Note 4:** The definitions of $D_1u_t(x,y)$, $D_2u_t(x,y)$ in M2 imply that: (i) For all y∈X, $D_1u_t(0,y)$ is the right-hand partial derivative with respect to the first variable at (0, y) and $D_1u_t(b,y)$ is the left-hand partial derivative with respect to the first variable at (b, y); (ii) For all x∈X, $D_2u_t(x,0)$ is the right-hand partial derivative with respect to the second variable at (x, 0) and $D_2u_t(x,b)$ is the left-hand partial derivative with respect to the second variable at (x, b).

For x∈X, let $\mathcal{F}(x) = \{<x_t|t\in\mathbb{N}^0>|\ (x_t, x_{t+1})\in\Omega_t, t\in\mathbb{N}^0, x_0 = x\}$.

We will (whenever necessary) refer to an infinite sequence $<x_t|t\in\mathbb{N}^0>\in \mathcal{F}(x)$ as a **trajectory starting at (from)** x.

Clearly, $\mathcal{F}(x)$ is a convex set for all x∈X.

Let $\mathcal{S}(x) = \{<x_t|t\in\mathbb{N}^0>\in \mathcal{F}(x)|\ \sum_{t=0}^{\infty} u_t(x_t, x_{t+1}) \geq \sum_{t=0}^{\infty} u_t(y_t, y_{t+1})$ for all $<y_t|t\in\mathbb{N}^0>\in \mathcal{F}(x)\}$.

Since, for all $t\in\mathbb{N}^0$, $u_t$:X×X→ ℝ is concave, it must be the case that $\mathcal{S}(x)$ is a convex set for all x∈X.

**Proposition 2.1:** Suppose M1, M2 and M3 are satisfied. Let $<x_t|t\in\mathbb{N}^0>\in \mathcal{F}(x)$. Then, $<x_t|t\in\mathbb{N}^0>\in\mathcal{S}(x)$ <u>if and only if</u> for all $<y_t|t\in\mathbb{N}^0>\in\mathcal{F}(x)$ it is the case that $D_1u_0(x_0, x_1)(x_0 - y_0) + \sum_{t=1}^{\infty}[D_2u_{t-1}(x_{t-1}, x_t) + D_1u_t(x_t, x_{t+1})](x_t - y_t) = \sum_{t=1}^{\infty}[D_2u_{t-1}(x_{t-1}, x_t) + D_1u_t(x_t, x_{t+1})](x_t - y_t) \geq 0$.

**Proof:** Let $<x_t|t\in\mathbb{N}^0>\in\mathcal{F}(x)$.

Note that for $\langle y_t | t \in \mathbb{N}^0 \rangle \in \mathcal{F}(x)$, $\sum_{t=0}^{\infty}[D_1u_t(x_t, x_{t+1})(x_t - y_t) + D_2u_t(x_t, x_{t+1})(x_{t+1} - y_{t+1})] = D_1u_0(x_0, x_1)(x_0 - y_0) + \sum_{t=1}^{\infty}[D_2u_{t-1}(x_{t-1}, x_t) + D_1u_t(x_t, x_{t+1})](x_t - y_t) = \sum_{t=1}^{\infty}[D_2u_{t-1}(x_{t-1}, x_t) + D_1u_t(x_t, x_{t+1})](x_t - y_t)$, since $y_0 = x_0 = x$.

Thus, for all $\langle x_t | t \in \mathbb{N}^0 \rangle$, $\langle y_t | t \in \mathbb{N}^0 \rangle \in \mathcal{F}(x)$, $D_1u_0(x_0, x_1)(x_0 - y_0) + \sum_{t=1}^{\infty}[D_2u_{t-1}(x_{t-1}, x_t) + D_1u_t(x_t, x_{t+1})](x_t - y_t) = \sum_{t=1}^{\infty}[D_2u_{t-1}(x_{t-1}, x_t) + D_1u_t(x_t, x_{t+1})](x_t - y_t) \geq 0$ <u>if and only if</u> $\sum_{t=0}^{\infty}[D_1u_t(x_t, x_{t+1})(x_t - y_t) + D_2u_t(x_t, x_{t+1})(x_{t+1} - y_{t+1})] \geq 0$.

Since for all $t \in \mathbb{N}^0$, $u_t: X \times X \to \mathbb{R}$ satisfies M2 and M3, for all $\langle y_t | t \in \mathbb{N}^0 \rangle \in \mathcal{F}(x)$ it must be the case that $\sum_{t=0}^{\infty} u_t(y_t, y_{t+1}) \leq \sum_{t=0}^{\infty} u_t(x_t, x_{t+1}) + \sum_{t=0}^{\infty}[D_1u_t(x_t, x_{t+1})(y_t - x_t) + D_2u_t(x_t, x_{t+1})(y_{t+1} - x_{t+1})]$.

If $\sum_{t=0}^{\infty}[D_1u_t(x_t, x_{t+1})(x_t - y_t) + D_2u_t(x_t, x_{t+1})(x_{t+1} - y_{t+1})] \geq 0$ for all $\langle y_t | t \in \mathbb{N}^0 \rangle \in \mathcal{F}(x)$, then $\sum_{t=0}^{\infty}[D_1u_t(x_t, x_{t+1})(y_t - x_t) + D_2u_t(x_t, x_{t+1})(y_{t+1} - x_{t+1})] \leq 0$ for all $\langle y_t | t \in \mathbb{N}^0 \rangle \in \mathcal{F}(x)$.

Thus, $\sum_{t=0}^{\infty} u_t(y_t, y_{t+1}) \leq \sum_{t=0}^{\infty} u_t(x_t, x_{t+1}) + \sum_{t=0}^{\infty}[D_1u_t(x_t, x_{t+1})(y_t - x_t) + D_2u_t(x_t, x_{t+1})(y_{t+1} - x_{t+1})]$ for all $\langle y_t | t \in \mathbb{N}^0 \rangle \in \mathcal{F}(x)$ implies $\sum_{t=0}^{\infty} u_t(x_t, x_{t+1}) \geq \sum_{t=0}^{\infty} u_t(y_t, y_{t+1})$ for all $\langle y_t | t \in \mathbb{N}^0 \rangle \in \mathcal{F}(x)$.

Thus, $\langle x_t | t \in \mathbb{N}^0 \rangle \in \mathcal{S}(x)$.

Now suppose, $\langle x_t | t \in \mathbb{N}^0 \rangle \in \mathcal{S}(x)$ and towards a contradiction suppose that for some $\langle y_t | t \in \mathbb{N} \rangle \in \mathcal{F}(x)$ it is the case that $\sum_{t=0}^{\infty}[D_1u_t(x_t, x_{t+1})(x_t - y_t) + D_2u_t(x_t, x_{t+1})(x_{t+1} - y_{t+1})] < 0$.

Clearly, $\langle x_t + \theta(y_t - x_t) | t \in \mathbb{N}^0 \rangle \in \mathcal{F}(x)$ for all $\theta \in [0, 1]$.

Consider the function $g: [0, 1] \to \mathbb{R}$ such that for $\theta \in [0, 1]$, $g(\theta) = \sum_{t=0}^{\infty} u_t(x_t + \theta(y_t - x_t), x_{t+1} + \theta(y_{t+1} - x_{t+1}))$.

$g(0) = \sum_{t=0}^{\infty} u_t(x_t, x_{t+1})$ and for all $\theta \in [0,1]$, $\frac{dg(\theta)}{d\theta} = \lim_{\alpha \to 0, (x_t + (\theta+\alpha)(y_t - x_t), x_{t+1} + (\theta+\alpha)(y_{t+1} - x_{t+1})) \in X \times X} \frac{g(\theta+\alpha) - g(\theta)}{\alpha}$.

By M3, $\frac{dg(\theta)}{d\theta} = \sum_{t=0}^{\infty}[D_1u_t(x_t + \theta(y_t - x_t), x_{t+1} + \theta(y_{t+1} - x_{t+1}))(y_t - x_t) + D_2u_t(x_t + \theta(y_t - x_t), x_{t+1} + \theta(y_{t+1} - x_{t+1}))(y_{t+1} - x_{t+1})]$.

Thus, g is continuously differentiable on [0,1].

Further, $\frac{dg(0)}{d\theta} = \sum_{t=0}^{\infty}[D_1u_t(x_t, x_{t+1})(y_t - x_t) + D_2u_t(x_t, x_{t+1})(y_{t+1} - x_{t+1})] > 0$.

By the continuous differentiability of g on [0,1], there exists $\alpha > 0$, such that for all $\theta \in [0, \alpha)$, $\frac{dg(\theta)}{d\theta} > 0$.

By the mean value theorem, there exists $\theta \in (0, \frac{\alpha}{2})$ such that $g(\frac{\alpha}{2}) - g(0) = \frac{\alpha}{2}\frac{dg(\theta)}{d\theta} > 0$.

Thus, $\sum_{t=0}^{\infty} u_t(x_t + \frac{\alpha}{2}(y_t - x_t), x_{t+1} + \frac{\alpha}{2}(y_{t+1} - x_{t+1})) = g(\frac{\alpha}{2}) > g(0) = \sum_{t=0}^{\infty} u_t(x_t, x_{t+1})$.

This contradicts our assumption that $\langle x_t | t \in \mathbb{N}^0 \rangle \in \mathcal{S}(x)$.

Hence, it must be the case that $\sum_{t=0}^{\infty}[D_1u_t(x_t,x_{t+1})(x_t - y_t) + D_2u_t(x_t,x_{t+1})\,(\,x_{t+1} - y_{t+1})] \geq 0$ for all $<y_t|t\in\mathbb{N}^0>\in\mathcal{F}(x)$.

Thus, $D_1u_0(x_0,x_1)(x_0 - y_0)+ \sum_{t=1}^{\infty}[D_2u_{t-1}(x_{t-1},x_t) + D_1u_t(x_t,x_{t+1})]\,(\,x_t - y_t) = \sum_{t=1}^{\infty}[D_2u_{t-1}(x_{t-1},x_t) + D_1u_t(x_t,x_{t+1})]\,(\,x_t - y_t) \geq 0$ for all $<y_t|t\in\mathbb{N}^0>\in\mathcal{F}(x)$. Q.E.D.

Proposition 2.1 provides the basis for a model "midway" between the reduced form model in Mitra (2000) and the "general dynamic optimization problem" we have considered so far. Note that for the initial results in Mitra (2000), differentiability is not assumed.

**3. Linear Dynamic Optimization:**

From this section onwards, we assume the following:

(M4) (i) For all $t\in\mathbb{N}^0$, $\Omega_t$ is a non-empty and closed subset of X×X; (ii) For all $t\in\mathbb{N}^0$ and x∈X, there exists y∈X (possibly depending on x) such that (x, y)$\in\Omega_t$. (iii) There exists an infinite sequence $<p^{(t)}|\, t\in\mathbb{N}^0>$ satisfying $\sum_{t=0}^{\infty}|p^{(t)}| < +\infty$ for such that for all x, y∈X and $t\in\mathbb{N}^0$, $u_t(x,y) = p^{(t)}x$.

**Note 5:** (iii) of M4 **<u>implies</u>** $\lim_{t\to\infty}|p^{(t)}| = 0$ and hence for all sequence $<x_t|t\in\mathbb{N}^0>$ with $x_t\in X$ for all $t\in\mathbb{N}^0$, it must be the case that $\lim_{t\to\infty}|p^{(t)}x_t| = 0$. Further, we <u>have not</u> assumed any convexity property for any $\Omega_t$.

Under assumption M4, we shall be concerned with the following problem:

Given x∈X: Maximize $\sum_{t=0}^{\infty}p^{(t)}x_t$, subject to the infinite sequence $<x_t|t\in\mathbb{N}^0>$ satisfying the constraints: $(x_t, x_{t+1})\in\Omega_t$, $t\in\mathbb{N}^0$, $x_0 = x$.

We shall refer to this problem (<u>satisfying M4</u>) as the **linear dynamic optimization (LDO) problem** and represent it as $<(p^{(t)}, \Omega_t)|\, t\in\mathbb{N}^0>$.

Kamihigashi and Roy (2006) consider dynamic optimization problems where for some real valued function f on X and δ∈(0, 1), the objective function is of the form $f(x_0) + \sum_{t=1}^{\infty}[f(x_t) - \delta^{t-1}x_t]$ for all sequences $<x_t|t\in\mathbb{N}^0>$ in X. Thus, they are not concerned with LDO problems in our sense, unless f is a linear function on X. If f is a linear function on X with slope μ, then the objective function would be linear with $p^{(0)} = \mu$ and $p^{(t)} = \mu - \delta^{t-1}$ for all $t\in\mathbb{N}$. In fact, we can generalize the linear case somewhat by introducing a sequence of real numbers $<\mu_t|\, t\in\mathbb{N}^0>$ such that $p^{(0)} = \mu_0$, $p^{(t)} = \mu_t - \delta^{t-1}$ for all $t\in\mathbb{N}$ and $\Omega_t = \{(x_t, x_{t+1})\in X\times X|\, x_{t+1} \leq \mu_t x_t\}$. In this case we arrive at a version of the model in Kamihigashi and Roy (2006) that is an LDO problem, since $\Omega_t$ is a non-empty and closed subset of X×X. The reason, why we felt it necessary to mention Kamihigashi and Roy (2006) is they refer to their objective function as linear, which is indeed the case, in the sense of "optimal growth" where for each $t\in\mathbb{N}^0$ an auxiliary variable $c_t = f(x_t) - x_{t+1}$ is introduced, whence their objective function can be reformulated as $\sum_{t=0}^{\infty}\delta^t c_t$.

**Note 6:** The kind of linearity discussed above can be used to formulate a "**linear model of optimal wealth accumulation**". Let $<u_t|t\in\mathbb{N}^0>$ be a sequence of non-negative real numbers that give the instantaneous average utility of monetary expenditure in each time-period and let $<\mu_t | t \in \mathbb{N}^0>$ be a sequence of non-negative real numbers that give the average rate of transformation of inherited wealth in each time-period.

The optimal wealth accumulation problem faced by the decision maker is the following:

Given x∈X: Maximize $\sum_{t=0}^{\infty} u_t(\mu_t x_t - x_{t+1})$, subject to $x_{t+1} \leq \mu_t x_t$, $x_t \in [0, b]$ for all $t \in \mathbb{N}^0$, $x_0 = x$.

For all $t \in \mathbb{N}^0$, expenditure during time-period t is given by $\mu_t x_t - x_{t+1}$.

The above problem reduces to the following LDO problem:

Given x∈X: Maximize $\sum_{t=0}^{\infty} p^{(t)} x_t$, subject to $x_{t+1} \leq \mu_t x_t$, $x_t \in [0, b]$ for all $t \in \mathbb{N}^0$, $x_0 = x$, where $p^{(0)} = u_0 \mu_0$ and $p^{(t)} = u_t \mu_t - u_{t-1}$ for all $t \in \mathbb{N}$.

In the case of this linear model of optimal wealth accumulation, for all $t \in \mathbb{N}^0$, $\Omega_t = \{(x, y) | x \in [0, b], y \in [0, \min\{\mu_t x, b\}]$.

**Note 7:** Another possible interpretation of the sequence $<p^{(t)} | t \in \mathbb{N}^0>$ could be associated with a sequence $<\pi^{(t)} | t \in \mathbb{N}^0>$ in the closed interval [0, 1] such that $\pi^{(0)} = 1$ and $\pi^{(T)} = 0$ for any $T \in \mathbb{N}$ implies $p^{(t)} = 0$ for all $t \geq T$. For $t \in \mathbb{N}$, $\pi^{(t)}$ denotes the probability that the optimization problem will "carry on" from period t -1 to period t. If instantaneous satisfaction in period 't' is derived from $(x_t - x_{t+1})$ and if the average instantaneous satisfaction in each period is given by the sequence $<v^{(t)} | t \in \mathbb{N}^0>$, then with $p^{(0)} = v^{(0)}$ and (i) $p^{(t)} = \pi^{(t)} v^{(t)} - \pi^{(t-1)} v^{(t-1)}$ for all $t \leq T = \max\{t | \pi^{(\tau)} > 0$ for all $\tau \leq t\}$ provided the maximum exists (and hence $T \in \mathbb{N}$), (ii) $p^{(t)} = \pi^{(t)} v^{(t)} - \pi^{(t-1)} v^{(t-1)}$ for all $t > 0$ otherwise, the LDO problem would imply "constrained" maximization of expected satisfaction.

**Note 8:** Based on the discussion in section 5.1 of Sorger (2015), we may consider a type of linear dynamic optimization problem which requires the existence of an infinite sequence $<(p_t^{(t)}, p_t^{(t+1)}) | t \in \mathbb{N}^0>$ of "pairs of real numbers" satisfying $\sum_{t=0}^{\infty} |p_t^{(t)}| < +\infty$ and $\sum_{t=0}^{\infty} |p_t^{(t+1)}| < +\infty$, such that $p^{(0)} = p_0^{(0)}$ and $p^{(t)} = p_t^{(t)} + p_{t-1}^{(t)}$ for all $t \in \mathbb{N}$.

We shall refer to a problem in the sub-class of optimization problems defined in note 8, as an **interlinked linear dynamic optimization (ILDO) problem** and represent it as $<((p_t^{(t)}, p_t^{(t+1)}), \Omega_t) | t \in \mathbb{N}^0>$.

We will refer to the LDO problem $<(p^{(t)}, \Omega_t) | t \in \mathbb{N}^0>$ such that $p^{(0)} = p_0^{(0)}$ and for all $t \in \mathbb{N}$, $p^{(t)} = p_t^{(t)} + p_{t-1}^{(t)}$ as the **canonical linear dynamic optimization (canonical LDO) problem associated with** $<((p_t^{(t)}, p_t^{(t+1)}), \Omega_t) | t \in \mathbb{N}^0>$.

**Note 9:** We will refer to an LDO problem that satisfies $\Omega_t = \{(x_t, x_{t+1}) \in X \times X | x_{t+1} \in [0, x_t]\}$ for all $t \in \mathbb{N}^0$ as a **linear cake-eating problem**.

Thus, for a linear cake-cutting problem $\Omega_t = \{(z, y) | y \in [0, z], z \in [0, b]\}$ for all $t \in \mathbb{N}^0$.

**For what follows we assume that $<((p^{(t)}, \Omega_t) | t \in \mathbb{N}^0>$ is a given LDO problem. As and when necessary, we will impose additional assumptions on this LDO problem.**

For all T∈$\mathbb{N}^0$, and y∈X, let $\mathcal{F}^T(y) = \{\langle x_t | t \geq T\rangle \mid (x_t, x_{t+1}) \in \Omega_t$ for all $t \geq T$ and $x_T = y\}$, and let $\mathcal{S}^T(y) = \underset{\langle x_t | t \geq T\rangle \in \mathcal{F}^T(y)}{\text{argmax}} \sum_{t=T}^{\infty} p^{(t)} x_t$.

Clearly, $\mathcal{F}^0(y) = \mathcal{F}(y)$ and $\mathcal{S}^0(y) = \mathcal{S}(y)$ for all y∈X.

For T∈$\mathbb{N}^0$ and y∈X, $\langle x_t | t \geq T\rangle \in \mathcal{F}^T(y)$ may be referred to as a **trajectory starting at (from)** y **at time-period** T.

It is easy to see that for all T∈$\mathbb{N}^0$ and y∈X, $\mathcal{F}^T(y)$ is non-empty.

**Remark:** An interesting property satisfied by linear models of "optimal wealth accumulation" is "**monotonicity of trajectories property**": For any T∈$\mathbb{N}^0$, and x, y∈X with y > x, [$\langle x_t | t \geq T\rangle \in \mathcal{F}^T(x)$, and $\langle y_t | t \geq T\rangle$ satisfies $y_T = y$, $y_t = x_t$ for all t > T] implies [$\langle y_t | t \geq T\rangle \in \mathcal{F}^T(y)$].

**4. Existence of solution, "Euler type" condition, competitive condition, and transversality condition:**

The proof of the following proposition is analogous to the proof of proposition 2.1 in Mitra (2000).

**Proposition 4.1:** $\mathcal{S}^T(x) \neq \phi$ for all x∈X.

**Proof:** Let (x, T)∈X×$\mathbb{N}^0$.

By (iii) of M4 and the assumption that $\Omega_t \subset X \times X = [0, b] \times [0, b]$ for all t∈$\mathbb{N}^0$ it must be the case that for all sequences $\langle x_t | t \in \mathbb{N}^0\rangle$ with $(x_t, x_{t+1}) \in X \times X$ for all t∈$\mathbb{N}^0$: $0 \leq \sum_{t=0}^{\infty} |p^{(t)} x_t| \leq b\sum_{t=0}^{\infty} |p^{(t)}| < +\infty$.

Let $M = b \sum_{t=0}^{\infty} |p^{(t)}| < +\infty$.

Let $\langle y_t | t \geq T\rangle \in \mathcal{F}^T(x)$. Clearly, $-M \leq \sum_{t=0}^{\infty} p^{(t)} y_t \leq M$.

Thus, $S(x) = \sup \{\sum_{t=T}^{\infty} p^{(t)} y_t \mid \langle y_t | t \geq T\rangle \in \mathcal{F}^T(x)\} < +\infty$.

Thus, there is a sequence of infinite sequences $\langle\langle x_t^{(n)} | t \geq T\rangle \mid n \in \mathbb{N}\rangle$ in $\mathcal{F}^T(x)$ such that for all n∈$\mathbb{N}$, $S(x) - \frac{1}{n} < \sum_{t=T}^{\infty} p^{(t)} x_t^{(n)}$.

Since, $\langle x_{T+1}^{(n)} \mid n \in \mathbb{N}\rangle$ is a sequence in the closed and bounded set X, it has a convergent subsequence $\langle x_{T+1}^{N_1(n)} | n \in \mathbb{N}\rangle$ converging to $x_{T+1}^0 \in X$.

Further, $(x, x_{T+1}^{N_1(n)}) \in \Omega_T$ for all n∈$\mathbb{N}$ implies $(x, x_{T+1}^0) \in \Omega_T$ since $\Omega_T$ is closed.

Consider the sequence $\langle x_{T+2}^{N_1(n)} | n \in \mathbb{N}\rangle$.

By an argument similar to the one in the previous step, it has a convergent subsequence $\langle x_{T+2}^{N_2(n)} \mid n \in \mathbb{N}\rangle$ converging to $x_{T+2}^0 \in X$.

Since, $\langle x_{T+1}^{N_2(n)} \mid n \in \mathbb{N}\rangle$ is a subsequence of the convergent subsequence $\langle x_{T+1}^{N_1(n)} \mid n \in \mathbb{N}\rangle$, it must be the case that $\langle x_{T+1}^{N_2(n)} \mid n \in \mathbb{N}\rangle$ converges to $x_{T+1}^0$.

The sequence $<(x_{T+1}^{N_2(n)}, x_{T+2}^{N_2(n)})|\ n\in\mathbb{N}>$ is in $\Omega_{T+1}$ and converges to $(x_{T+1}^0, x_{T+2}^0)$. Since $\Omega_{T+1}$ is closed $(x_{T+1}^0, x_{T+2}^0)\in \Omega_{T+1}$.

Having obtained convergent subsequences $<(x_{T+\tau-1}^{N_\tau(n)}, x_{T+\tau}^{N_\tau(n)})|\ n\in\mathbb{N}>\in\Omega_{T+\tau-1}$ for all $\tau = 1, \ldots, t$ for some $t \geq 1$, let $<(x_{T+t}^{N_{t+1}(n)}, x_{T+t}^{N_{t+1}(n)}|\ n\in\mathbb{N}>$ be a convergent subsequence of the sequence $<(x_{T+t}^{N_t(n)}, x_{T+t+1}^{N_t(n)}|\ n\in\mathbb{N}>$ in the closed and bounded set X×X, converging to $(x_{T+t}^0, x_{T+t+1}^0)$ $\in$X×X×X. Since, $<(x_{T+t}^{N_t(n)}, x_{T+t+1}^{N_t(n)}|\ n\in\mathbb{N}>\in\Omega_{T+t}$ and $\Omega_{T+t}$ is a closed subset of X×X, it must be the case that $(x_{T+t}^0, x_{T+t+1}^0)\in \Omega_{T+t}$.

Consider the sequence $<x_{T+t}^0|t\in\mathbb{N}^0>$. Since, $(x_{T+t}^0, x_{T+t+1}^0)\in\Omega_t$ for all $t\in\mathbb{N}^0$, where $x_T^0$ = x, it must be the case that $<x_{T+t}^0|t\in\mathbb{N}^0>\in\mathcal{F}(x)$.

Thus, $\sum_{t=0}^{\infty} p^{(T+t)}x_{T+t}^0 \leq S(x)$.

We wish to show that $\sum_{t=0}^{\infty} p^{(T+t)}x_{T+t}^0 = S(x)$.

Since, for all $t\in\mathbb{N}^0$, $<x_{T+\tau}^{N_\tau(n)}|\ n\in\mathbb{N}>$ converges to $x_{T+\tau}^0$ and $(x_{T+\tau}^{N_\tau(n)}, x_{T+\tau+1}^{N_\tau(n)})\in \Omega_{T+\tau}$, for all $\tau \leq t$ and $n\in\mathbb{N}$, ***given*** $\varepsilon > 0$, for all $t\in\mathbb{N}$, there exists a sequence $<n_t|t\in \mathbb{N}>$ in $\mathbb{N}$ satisfying $n_{t+1} > n_t$ for all $t\in\mathbb{N}$ such that for all $n\in\mathbb{N}$ with $n \geq n_t$ and $\tau \leq t$, $(x_{T+\tau}^{N_t(n)}, x_{T+\tau+1}^{N_t(n)})\in \Omega_{T+\tau}$ and $|p^{(T+\tau)}x_{T+\tau}^0 - p^{(T+\tau)}x_{T+\tau}^{N_t(n)}| < \frac{\varepsilon}{8}(\frac{1}{2})^t$, i.e, $p^{(T+\tau)}x_{T+\tau}^{N_t(n)} + \frac{\varepsilon}{8}(\frac{1}{2})^t > p^{(T+\tau)}x_{T+\tau}^0 > p^{(T+\tau)}x_{T+\tau}^{N_t(n)} - \frac{\varepsilon}{8}(\frac{1}{2})^t$.

Note that $\sum_{\tau=1}^{t}\frac{\varepsilon}{8}(\frac{1}{2})^t \leq \sum_{\tau=1}^{t}\frac{\varepsilon}{8}(\frac{1}{2})^\tau = \frac{\varepsilon}{16}\sum_{\tau=0}^{t-1}(\frac{1}{2})^\tau = \frac{\varepsilon}{8}(1- (\frac{1}{2})^t)$ and $x_T^{N_t(n)}$ = x for all $(t, n)\in\mathbb{N}\times\mathbb{N}$.

Thus, for all $t\in\mathbb{N}$, $\sum_{\tau=0}^{t} p^{(T+\tau)}x_{T+\tau}^0 > \sum_{\tau=0}^{t} p^{(T+\tau)}x_{T+\tau}^{N_t(n)} - \frac{\varepsilon}{8}(1- (\frac{1}{2})^t)$ for all $n \geq n_t$.

For all $t\in\mathbb{N}$: $\sum_{\tau=0}^{t} p^{(T+\tau)}x_{T+\tau}^{N_t(n)} > S(x) - \frac{1}{N_t(n)} - \sum_{\tau=t+1}^{\infty} p^{(T+\tau)}x_{T+\tau}^{N_t(n)}$ for all $n \geq n_t$.

Thus, for all $t\in\mathbb{N}$, $\sum_{\tau=0}^{t} p^{(T+\tau)}x_{T+\tau}^0 > \sum_{\tau=0}^{t} p^{(T+\tau)}x_{T+\tau}^{N_t(n)} - \frac{\varepsilon}{8}(1- (\frac{1}{2})^t) > S(x) - \frac{1}{N_t(n)} - \sum_{\tau=t+1}^{\infty} p^{(T+\tau)}x_{T+\tau}^{N_t(n)} - \frac{\varepsilon}{8}(1- (\frac{1}{2})^t)$ for all $n \geq n_t$.

Thus, for all $t\in\mathbb{N}$, $\sum_{\tau=0}^{t} p^{(T+\tau)}x_{T+\tau}^0 + \sum_{\tau=t+1}^{\infty} p^{(T+\tau)}x_{T+\tau}^{N_t(n)} > S(x) - \frac{1}{N_t(n)} - \frac{\varepsilon}{8}(1- (\frac{1}{2})^t)$ for all $n \geq n_t$.

By hypothesis, $\sum_{t=0}^{\infty} |p^{(t)}| < +\infty$ and for all $(t, n)\in \mathbb{N}^0\times\mathbb{N}$, $x_{T+t}^{(n)}\in[0, b]$.

Thus, for all $t\in \mathbb{N}$ and $n \geq n_t$, $\sum_{\tau=t+1}^{\infty} |p^{(T+\tau)}x_{T+\tau}^{N_t(n)}| = \sum_{\tau=t+1}^{\infty} |p^{(T+\tau)}|x_{T+\tau}^{N_t(n)}$.

For all $t\in\mathbb{N}$ and $n\in\mathbb{N}$, $x_t^{(n)}\in[0, b]$ <u>implies</u> for all $t\in\mathbb{N}$ and $\tau \geq t + 1$, $\sup_{n\geq n_t} x_{T+\tau}^{N_t(n)}$ exists, $\sup_{n\geq n_t} x_{T+\tau}^{N_t(n)}\in[0, b]$ and $\sup_{n\geq n_t} x_{T+\tau}^{N_t(n)} \geq x_{T+\tau}^{N_t(n)}$ for all $n \geq n_t$.

Thus, for all $t\in\mathbb{N}$ and $n \geq n_t$, $\sum_{\tau=t+1}^{\infty}|p^{(T+\tau)}|[\sup_{n\geq n_t} x_{T+\tau}^{N_t(n)}] \geq \sum_{\tau=t+1}^{\infty} |p^{(T+\tau)}|x_{T+\tau}^{N_t(n)}$.

Further, for all $t\in\mathbb{N}$ and $n \geq n_t$, $\sum_{\tau=t+1}^{\infty} |p^{(T+\tau)}x_{T+\tau}^{N_t(n)}| \geq \sum_{\tau=t+1}^{\infty} p^{(T+\tau)}x_{T+\tau}^{N_t(n)}$.

Thus, for all t∈ℕ, $\sum_{\tau=0}^{t} p^{(T+\tau)} x_{T+\tau}^{0} + \sum_{\tau=t+1}^{\infty} |p^{(T+\tau)}|[\sup_{n \geq n_t} x_{T+\tau}^{N_t(n)}] \geq \sum_{\tau=1}^{t} p^{(T+\tau)} x_{T+\tau}^{0} + \sum_{\tau=t+1}^{\infty} p^{(T+\tau)} x_{T+\tau}^{N_t(n)} >$ S(x) - $\frac{1}{N_t(n)}$ - $\frac{\varepsilon}{8}(1-(\frac{1}{2})^t)$ for all n ≥ $n_t$.

Thus, for all t∈ℕ, $\sum_{\tau=0}^{t} p^{(T+\tau)} x_{T+\tau}^{0} + \sum_{\tau=t+1}^{\infty} |p^{(T+\tau)}|[\sup_{n \geq n_t} x_{T+\tau}^{N_t(n)}] >$ S(x) - $\frac{1}{N_t(n)}$ - $\frac{\varepsilon}{8}(1-(\frac{1}{2})^t)$ for all n ≥ $n_t$.

Thus, for all t∈ℕ, $\sum_{\tau=0}^{t} p^{(T+\tau)} x_{T+\tau}^{0} + \sum_{\tau=t+1}^{\infty} |p^{(T+\tau)}|[\sup_{n \geq n_t} x_{T+\tau}^{N_t(n)}] \geq$ S(x) - $\lim_{n\to\infty} \frac{1}{N_t(n)}$ - $\frac{\varepsilon}{8}(1-(\frac{1}{2})^t)$ = S(x) - $\frac{\varepsilon}{8}(1-(\frac{1}{2})^t)$, since $\lim_{n\to\infty} \frac{1}{N_t(n)} = 0$.

Thus, for all t∈ℕ, $\sum_{\tau=0}^{t} p^{(T+\tau)} x_{T+\tau}^{0} + \sum_{\tau=t+1}^{\infty} |p^{(T+\tau)}|[\sup_{n \geq n_t} x_{T+\tau}^{N_t(n)}] \geq$ S(x) - $\frac{\varepsilon}{8}(1-(\frac{1}{2})^t)$.

Since $\lim_{t\to\infty} \frac{\varepsilon}{8}(1-(\frac{1}{2})^t) = \frac{\varepsilon}{8}$ and $\lim_{t\to\infty} \sum_{\tau=0}^{t} p^{(T+\tau)} x_{T+\tau}^{0} = \sum_{\tau=0}^{\infty} p^{(T+\tau)} x_{T+\tau}^{0}$, we get $\sum_{\tau=0}^{\infty} p^{(T+\tau)} x_{T+\tau}^{0} + \lim_{t\to\infty} \sum_{\tau=t+1}^{\infty} |p^{(T+\tau)}|[\sup_{n \geq n_t} x_{T+\tau}^{N_t(n)}] \geq$ S(x) - $\frac{\varepsilon}{8}$.

Since, for all t∈ℕ, $\sum_{\tau=t+1}^{\infty} |p^{(T+\tau)}| < +\infty$ and $[\sup_{n \geq n_t} x_{T+\tau}^{N_t(n)}] \in [0, b]$, it must be the case that $\lim_{t\to\infty} \sum_{\tau=t+1}^{\infty} |p^{(T+\tau)}|[\sup_{n \geq n_t} x_{T+\tau}^{N_t(n)}] = 0$.

Thus, $\sum_{t=T}^{\infty} p^{(t)} x_t^0 \geq$ S(x) - $\frac{\varepsilon}{8}$.

Since the above holds for all ε > 0, we get $\sum_{t=T}^{\infty} p^{(t)} x_t^0 \geq$ S(x).

Thus, $\sum_{t=T}^{\infty} p^{(t)} x_t^0$ = S(x). Q. E. D.

**Note 9:** For the iterative process outlined above we are unaware of any possibility of it "invariably" converging to a strictly increasing sequence $<N_\infty(n)|n\in\mathbb{N}>$ such that for all t∈$\mathbb{N}^0$, $(x_t^{N_\infty(n)}, x_{t+1}^{N_\infty(n)}) \in \Omega_t$ for all n∈ℕ and $<x_t^{N_\infty(n)}| n\in\mathbb{N}>$ converges to $x_t^0$. However, in appendix A.2 (entitled "A Diagonal Process of Cantor") of Blot and Hayek (2014), it seems to be indicated (provided our understanding is correct), that for all t∈$\mathbb{N}^0$, the subsequence $<x_{T+t}^{N_n(n)}| n\in\mathbb{N}>$ derived from the construction in the proof above, converges to $x_{T+t}^0$.

By proposition 4.1 we know that for all (x, T)∈X×$\mathbb{N}^0$, $\mathcal{S}^T$(x) ≠ ϕ, and hence for all T∈$\mathbb{N}^0$, there exists a function $V^T$: X →ℝ such that for all x∈X, $V^T(x) = \sum_{t=T}^{\infty} p^{(t)} x_t$ for all $<x_{T+t}|t\in\mathbb{N}^0>\in\mathcal{S}^T$(x). Clearly for all for all T∈$\mathbb{N}^0$, $V^T$ is well defined on X.

We may refer to $V^T$ as the **period-T optimal value function**

We will denote $V^0$ by V. V is said to be the **optimal value function**.

Thus, for all x∈X, V(x) = $p_0^{(0)} x + \sum_{t=1}^{\infty} p^{(t)} x_t$, where $<x_t|t\in\mathbb{N}^0>\in\mathcal{S}$(x)

**Proposition 4.2:** For x∈X, let $<x_t|t\in\mathbb{N}^0>\in\mathcal{S}(x)$.

(i) Suppose that for some $T \in \mathbb{N}$, $\{y \in \Omega_{T-1}(x_{T-1}) | y > x_T\} \neq \phi$. If given $\varepsilon > 0$, there exists $y \in \Omega_{T-1}(x_{T-1})$ satisfying $y > x_T$ and there exists $<y_{T+t}|t \in \mathbb{N}^0> \in \mathcal{F}^T(y)$ such that $|y_{T+t} - x_{T+t}| < \varepsilon|y - x_T|$ for all $t \in \mathbb{N}$], then $p^{(T)} \leq 0$.

(ii) Suppose that for some $T \in \mathbb{N}$, $\{y \in \Omega_{T-1}(x_{T-1}) | y < x_T\} \neq \phi$. If given $\varepsilon > 0$, there exists $y \in \Omega_{T-1}(x_{T-1})$ satisfying $y < x_T$ and there exists $<y_{T+t}|t \in \mathbb{N}^0> \in \mathcal{F}^T(y)$ such that $|y_{T+t} - x_{T+t}| < \varepsilon|y - x_T|$ for all $t \in \mathbb{N}$], then $p^{(T)} \geq 0$.

(iii) $\lim_{t\to\infty} |p^{(t)}|x_t = 0$.

**Proof:** (i) Under the conditions stated in (i) of the proposition and towards a contradiction suppose $p^{(T)} > 0$. Let $\varepsilon > 0$ satisfy $p^{(T)} > \varepsilon\sum_{t=0}^{\infty} | p^{(t)}|$. Let $y \in \Omega_{T-1}(x_{T-1})$ satisfy $y > x_T$ and $<y_{T+t}|t \in \mathbb{N}^0> \in \mathcal{F}^T(y)$ such that $|(y_{T+t} - x_{T+t})| < \varepsilon|y - x_T|$ for all $t \in \mathbb{N}$.

Let $<z_t|t \in \mathbb{N}^0> \in \mathcal{F}(x)$ such that $z_t = x_t$ for $t = 0, \ldots, T-1$, $z_t = y_t$ for all $t \geq T$.

Thus, $\sum_{t=0}^{\infty} p^{(t)}z_t = \sum_{t=0}^{T-1} p^{(t)}x_t + p^{(T)}y + \sum_{t=T+1}^{\infty} p^{(t)}y_t = \sum_{t=0}^{T-1} p^{(t)}x_t + p^{(T)}x_T + p^{(T)}(y - x_T) + \sum_{t=T+1}^{\infty} p^{(t)}y_t$.

If for $t > T$, $p^{(t)} \geq 0$, then $p^{(t)}y_t \geq p^{(t)}x_t - \varepsilon p^{(t)}|y - x_T|$.

If for $t > T$, $p^{(t)} < 0$, then $p^{(t)}y_t > p^{(t)}x_t + \varepsilon p^{(t)}|y - x_T|$.

Thus, $\sum_{t=0}^{\infty} p^{(t)}z_t > \sum_{t=0}^{T-1} p^{(t)}x_t + p^{(T)}x_T + p^{(T)}(y - x_T) + \sum_{t=T+1}^{\infty} p^{(t)}y_t \geq \sum_{t=0}^{T-1} p^{(t)}x_t + p^{(T)}x_T + p^{(T)}(y - x_T) + \sum_{\{t>T|p^{(t)}\geq 0\}} p^{(t)}x_t - \varepsilon|y - x_T|\sum_{\{t>T|p^{(t)}\geq 0\}} p^{(t)} + + \sum_{\{t>T|p^{(t)}< 0\}} p^{(t)}x_t + \varepsilon|y - x_T|\sum_{\{t>T|p^{(t)}< 0\}} p^{(t)} = \sum_{t=0}^{T-1} p^{(t)}x_t + p^{(T)}x_T + p^{(T)}(y - x_T) + \sum_{t=T+1}^{\infty} p^{(t)}x_t - \varepsilon|y - x_T|\sum_{t=T+1}^{\infty} |p^{(t)}| = \sum_{t=0}^{\infty} p^{(t)}x_t + p^{(T)}(y - x_T) - \varepsilon|y - x_T|\sum_{t=T+1}^{\infty} |p^{(t)}| \geq \sum_{t=0}^{\infty} p^{(t)}x_t + p^{(T)}(y - x_T) - \varepsilon|y - x_T|\sum_{t=0}^{\infty} |p^{(t)}| = \sum_{t=0}^{\infty} p^{(t)}x_t + |y - x_T|(p^{(T)} - \varepsilon \sum_{t=0}^{\infty} |p^{(t)}|) > \sum_{t=0}^{\infty} p^{(t)}x_t$, since we have assumed $|y - x_T| = y - x_T > 0$ and $p^{(T)} > \varepsilon \sum_{t=0}^{\infty} |p^{(t)}|$.

Since $<z_t|t \in \mathbb{N}^0> \in \mathcal{F}(x)$ the assumption that $<x_t|t \in \mathbb{N}^0> \in \mathcal{S}(x)$ stands contradicted.

Hence, it must be the case that $p^{(T)} \leq 0$.

(ii) Under the conditions stated in (ii) of the proposition and towards a contradiction suppose $p^{(T)} < 0$. Let $\varepsilon > 0$ satisfy $- p^{(T)} > \varepsilon\sum_{t=0}^{\infty} | p^{(t)}|$. Let $y \in \Omega_{T-1}(x_{T-1})$ satisfy $y < x_T$ and $<y_{T+t}|t \in \mathbb{N}^0> \in \mathcal{F}^T(y)$ such that $|(y_{T+t} - x_{T+t})| < \varepsilon|y - x_T|$ for all $t \in \mathbb{N}$.

Let $<z_t|t \in \mathbb{N}^0> \in \mathcal{F}(x)$ such that $y_t = x_t$ for $t = 0, \ldots, T-1$, $z_t = y_t$ for all $t \geq T$.

Thus, $\sum_{t=0}^{\infty} p^{(t)}z_t = \sum_{t=0}^{T-1} p^{(t)}x_t + p^{(T)}y + \sum_{t=T+1}^{\infty} p^{(t)}y_t = \sum_{t=0}^{T-1} p^{(t)}x_t + p^{(T)}x_T + p^{(T)}(y - x_T) + \sum_{t=T+1}^{\infty} p^{(t)}y_t$.

If for $t > T$, $p^{(t)} \geq 0$, then $p^{(t)}y_t \geq p^{(t)}x_t - \varepsilon p^{(t)}|y - x_T|$.

If for $t > T$, $p^{(t)} < 0$, then $p^{(t)}y_t > p^{(t)}x_t + \varepsilon p^{(t)}|y - x_T|$.

Thus, $\sum_{t=0}^{\infty} p^{(t)}z_t = \sum_{t=0}^{T-1} p^{(t)}x_t + p^{(T)}x_T + p^{(T)}(y - x_T) + \sum_{t=T+1}^{\infty} p^{(t)}y_t \geq \sum_{t=0}^{T-1} p^{(t)}x_t + p^{(T)}x_T + p^{(T)}(y - x_T) + \sum_{\{t>T|p^{(t)}\geq 0\}} p^{(t)}x_t - \varepsilon|y - x_T|\sum_{\{t>T|p^{(t)}\geq 0\}} p^{(t)} + \sum_{\{t>T|p^{(t)}< 0\}} p^{(t)}x_t + \varepsilon|y - x_T|\sum_{\{t>T|p^{(t)}< 0\}} p^{(t)} = \sum_{t=0}^{T-1} p^{(t)}x_t + p^{(T)}x_T + p^{(T)}(y - x_T) + \sum_{t=T+1}^{\infty} p^{(t)}x_t - \varepsilon|y -$

$x_T|\sum_{t=T+1}^{\infty}|p^{(t)}| = \sum_{t=0}^{\infty} p^{(t)}x_t + p^{(T)}(y - x_T) - \varepsilon|y - x_T|\sum_{t=T+1}^{\infty}|p^{(t)}| \geq \sum_{t=0}^{\infty} p^{(t)}x_t + p^{(T)}(y - x_T) - \varepsilon|y - x_T|\sum_{t=0}^{\infty}|p^{(t)}| = \sum_{t=0}^{\infty} p^{(t)}x_t - |y - x_T|(- p^{(T)} - \varepsilon \sum_{t=0}^{\infty}|p^{(t)}|) > \sum_{t=0}^{\infty} p^{(t)}x_t$, since we have assumed $y - x_T < 0$, so that $|y - x_T| = -(y-x_T)$ and $-p^{(T)} > \varepsilon \sum_{t=0}^{\infty}|p^{(t)}|$.

Since $<z_t|t\in\mathbb{N}^0>\in\mathcal{F}(x)$ the assumption that $<x_t|t\in\mathbb{N}^0>\in\mathcal{S}(x)$ stands contradicted.

Hence, it must be the case that $p^{(T)} \leq 0$.

(iii) $\sum_{t=0}^{\infty}|p^{(t)}| < +\infty$ implies $\lim_{t\to\infty} |p^{(t)}| = 0$.

Further $x_t \in [0, b]$ for all $t\in\mathbb{N}^0$.

Thus, $\lim_{t\to\infty} |p^{(t)}|x_t = 0$.    Q.E.D.

**Note 10:** That $\{y\in\Omega_{T-1}(x_{T-1})|\ y > x_T\} \neq \phi$ under certain mild conditions implies $p^{(T)} \leq 0$ & that $\{y\in\Omega_{T-1}(x_{T-1})|\ y < x_T\} \neq \phi$ under certain mild conditions implies $p^{(T)} \geq 0$, are similar to what is referred to in section 5.2 of Sorger (2015) as the Euler equation. We may refer to "$\{y\in\Omega_{T-1}(x_{T-1})|\ y > x_T\} \neq \phi$ under certain mild conditions implies $p^{(T)} \leq 0$ & $\{y\in\Omega_{T-1}(x_{T-1})|\ y < x_T\} \neq \phi$ under certain mild conditions implies $p^{(T)}$" as "Euler conditions." The condition $\lim_{T\to\infty} |p^{(T)}|x_T = 0$ resembles the several versions of "transversality condition" that are available in the literature on dynamic optimization.

The following result, analogous to proposition 4.1 in Mitra (2000), provides a "neat" sufficient condition for a trajectory starting at x to belong to $\mathcal{S}(x)$.

**Proposition 4.3:** If for $<x_t|t\in\mathbb{N}^0>\in\mathcal{F}(x)$ and there is a sequence $<q_t|t\in\mathbb{N}^0>$ of non-negative real numbers satisfying:

(i) For all $t\in\mathbb{N}^0$: $p^{(t)}x_t + q_{t+1}x_{t+1} - q_tx_t \geq p^{(t)}y + q_{t+1}z - q_ty$ for all $(y, z)\in\Omega_t$,

(ii) $\lim_{t\to\infty} q_tx_t = 0$,

then, $<x_t|t\in\mathbb{N}^0>\in\mathcal{S}(x)$.

**Proof:** Let $<y_t|t\in\mathbb{N}^0>\in\mathcal{F}(x)$.

Then, (i) implies that for all $t\in\mathbb{N}^0$: $p^{(t)}x_t + q_{t+1}x_{t+1} - q_tx_t \geq p^{(t)}y_t + q_{t+1}y_{t+1} - q_ty_t$ and thus for all $t\in\mathbb{N}^0$: $p^{(t)}x_t - p^{(t)}y_t \geq (q_{t+1}y_{t+1} - q_ty_t) - (q_{t+1}x_{t+1} - q_tx_t)$.

Hence, for all $T\in\mathbb{N}^0$: $\sum_{t=0}^{T} p^{(t)}x_t - \sum_{t=0}^{T} p^{(t)}y_t \geq \sum_{t=0}^{T}[q_{t+1}y_{t+1} - q_ty_t] - \sum_{t=0}^{T}[q_{t+1}x_{t+1} - q_tx_t] = q_{T+1}y_{T+1} - q_{T+1}x_{T+1} \geq - q_{T+1}x_{T+1}$, since $q_t \geq 0$ and $y_t \geq 0$ for all $t\in\mathbb{N}^0$.

Thus, $\sum_{t=0}^{\infty} p^{(t)}x_t - \sum_{t=0}^{\infty} p^{(t)}y_t = \lim_{T\to\infty}\sum_{t=0}^{T} p^{(t)}x_t - \lim_{T\to\infty}\sum_{t=0}^{T} p^{(t)}y_t \geq \lim_{T\to\infty} -q_{T+1}x_{T+1} = 0$, by (ii). Q.E.D.

**Note 11:** Conditions (i), is referred to as the "competitive condition" in Mitra (2000). The competitive condition *can also be written as* "for all $t\in\mathbb{N}^0$: $[V^t(x_t) - q_tx_t] - [V^t(x_{t+1}) - q_{t+1}x_{t+1}] \geq [V^t(y) - q_ty] - [V^{t+1}(z) - q_{t+1}z]$ for all $(y, z)\in\Omega_t$."

**5. Linear Dynamic Programming:**

The following propositions are similar to proposition 3.1 in Mitra (2000).

**Proposition 5.1:** Suppose $\Omega_T$ is convex for all $T \in \mathbb{N}^0$. Then, for all $T \in \mathbb{N}^0$, $V^T$ is concave and continuous on X.

**Proof:** Let $T \in \mathbb{N}^0$. For all (x, T), (y, T)$\in X \times \mathbb{N}^0$, let $<x_t | t \geq T> \in \mathcal{S}^T(x)$ and $<y_t | t \geq T> \in \mathcal{S}^T(y)$. Suppose $\theta \in (0, 1)$ and consider the infinite sequence $<\theta x_t + (1-\theta) y_t | t \geq T>$.

Since $(x_t, x_{t+1}) \in \Omega_t$ and $(y_t, y_{t+1}) \in \Omega_t$ for all $t \geq T$ and $\Omega_t$ is a convex set for all $t \geq T$, $(\theta x_t + (1-\theta) y_t, \theta x_{t+1} + (1-\theta) y_{t+1}) \in \Omega_t$ for all $t \geq T$.

Thus, $<\theta x_t + (1-\theta) y_t | t \geq T> \in \mathcal{F}^T(\theta x + (1-\theta) y)$.

Thus, $V^T(\theta x + (1-\theta)y) \geq p_T^{(T)}(\theta x + (1-\theta)y) + \sum_{t=T+1}^{\infty} p^{(t)}(\theta x_t + (1-\theta)y_t) = \theta[p_T^{(T)} x + \sum_{t=T+1}^{\infty} p^{(t)} x_t] + (1-\theta)[p_T^{(T)} y + \sum_{t=T+1}^{\infty} p^{(t)} y_t] = \theta V^T(x) + (1-\theta) V^T(y)$.

Thus, for all $T \in \mathbb{N}^0$, $V^T$ is a concave function on X.

Since $V^T$ is concave on X = [0, b], $V^T$ is continuous on (0, b) and both $\lim_{x \to 0} V(x)$ and $\lim_{x \to b} V(x)$ exists and belong to the interval [-M, M] where $M = b\sum_{t=0}^{\infty} |p^{(t)}| < +\infty$.

Further, it must be that $\lim_{x \to 0} V^T(x) \geq V^T(0)$ and $\lim_{x \to b} V^T(x) \geq V^T(b)$.

The results stated above follow from corollary 1 of proposition 1 and part (i) of proposition 3 in Lahiri (2025a).

Now let $x \in X$ and $<x^{(n)} | n \in \mathbb{N}>$ be a sequence in X converging to x.

For each $n \in \mathbb{N}$, let $<x_{T+t}^{(n)} | t \in \mathbb{N}^0> \in \mathcal{S}(x^{(n)})$.

Since $x_T^{(n)} = x^{(n)}$ for all $n \in \mathbb{N}$, $\lim_{n \to \infty} x_T^{(n)} = x$.

Since, $<x_{T+1}^{(n)} | n \in \mathbb{N}>$ is a sequence in the closed and bounded set X, it has a convergent subsequence $<x_{T+1}^{N_1(n)} | n \in \mathbb{N}>$ converging to $x_{T+1}^0 \in X$.

Further, $(x_T^{N_1(n)}, x_{T+1}^{N_1(n)}) \in \Omega_T$ for all $n \in \mathbb{N}$ implies $(x, x_{T+1}^0) \in \Omega_T$ since $\Omega_T$ is closed.

Having obtained convergent subsequences $<(x_{T+\tau-1}^{N_\tau(n)}, x_{T+\tau}^{N_\tau(n)}) | n \in \mathbb{N}> \in \Omega_{T+\tau-1}$ for all $\tau = 1, \ldots, t$ for some $t \geq 1$, let $<(x_{T+t}^{N_{t+1}(n)}, x_{T+t}^{N_{t+1}(n)}) | n \in \mathbb{N}>$ be a convergent subsequence of the sequence $<(x_{T+t}^{N_t(n)}, x_{T+t+1}^{N_t(n)}) | n \in \mathbb{N}>$ in the closed and bounded set $X \times X$, converging to $(x_{T+t}^0, x_{T+t+1}^0) \in X \times X$. Since, $<(x_{T+t}^{N_t(n)}, x_{T+t+1}^{N_t(n)}) | n \in \mathbb{N}> \in \Omega_{T+t}$ and $\Omega_{T+t}$ is a closed subset of $X \times X$, it must be the case that $(x_{T+t}^0, x_{T+t+1}^0) \in \Omega_{T+t}$.

Since, for all $t \in \mathbb{N}^0$, $<x_{T+\tau}^{N_t(n)} | n \in \mathbb{N}>$ converges to $x_{T+\tau}^0$ and $(x_{T+\tau}^{N_\tau(n)}, x_{T+\tau+1}^{N_\tau(n)}) \in \Omega_{T+\tau}$, for all $\tau \leq t$ and $n \in \mathbb{N}$, ***given*** $\varepsilon > 0$, for all $t \in \mathbb{N}$, there exists a sequence $<n_t | t \in \mathbb{N}>$ in $\mathbb{N}$ satisfying $n_{t+1} > n_t$ for all $t \in \mathbb{N}$ such that for all $n \in \mathbb{N}$ with $n \geq n_t$ and $\tau \leq t$, $(x_{T+\tau}^{N_t(n)}, x_{T+\tau+1}^{N_t(n)}) \in \Omega_{T+\tau}$ and $|p^{(T+\tau)} x_{T+\tau}^0 - p^{(T+\tau)} x_{T+\tau}^{N_t(n)}| < \frac{\varepsilon}{8}(\frac{1}{2})^t$, i.e, $p^{(T+\tau)} x_{T+\tau}^{N_t(n)} + \frac{\varepsilon}{8}(\frac{1}{2})^t > p^{(T+\tau)} x_{T+\tau}^0 > p^{(T+\tau)} x_{T+\tau}^{N_t(n)} - \frac{\varepsilon}{8}(\frac{1}{2})^t$.

Note that $\sum_{\tau=0}^{t}\frac{\varepsilon}{8}(\frac{1}{2})^{t} \le \sum_{\tau=0}^{t}\frac{\varepsilon}{8}(\frac{1}{2})^{\tau} = \frac{\varepsilon}{4}(1- (\frac{1}{2})^{t+1})$.

Thus, for all t∈ ℕ, $\sum_{\tau=0}^{t} p^{(T+\tau)}x_{T+\tau}^{0} > \sum_{\tau=0}^{t} p^{(T+\tau)}x_{T+\tau}^{N_\tau(n)} - \frac{\varepsilon}{4}(1- (\frac{1}{2})^{t+1})$ for all $n \ge n_t$.

Thus, for all t∈ℕ, $\sum_{\tau=0}^{t} p^{(T+\tau)}x_{T+\tau}^{0} > \sum_{\tau=0}^{t} p^{(T+\tau)}x_{T+\tau}^{N_t(n)} - \frac{\varepsilon}{4}(1- (\frac{1}{2})^{t+1}) = V(x^{N_t(n)}) - \sum_{t=\tau+1}^{\infty} p^{(T+t)}x_{T+t}^{N_\tau(n)} - \frac{\varepsilon}{4}(1- (\frac{1}{2})^{t+1})$ for all $n \ge n_t$.

Thus, for all t∈ℕ, $V^T(x) \ge \sum_{\tau=0}^{\infty} p^{(T+t)}x_{T+\tau}^{0} = \sum_{\tau=0}^{t} p^{(T+t)}x_{T+\tau}^{0} + \sum_{\tau=t+1}^{\infty} p^{(T+\tau)}x_{T+\tau}^{0} > V(x^{N_t(n)}) - \sum_{\tau=t+1}^{\infty} p^{(T+\tau)}x_{T+\tau}^{N_t(n)} - \frac{\varepsilon}{4}(1- (\frac{1}{2})^{t+1}) + \sum_{\tau=t+1}^{\infty} p^{(T+\tau)}x_{T+\tau}^{0}$ for all $n \ge n_\tau$.

Thus, for all t∈ℕ, $V^T(x) + \sum_{\tau=\tau+1}^{\infty} p^{(T+\tau)}x_{T+\tau}^{N_t(n)} > V^T(x^{N_t(n)}) - \frac{\varepsilon}{4}(1- (\frac{1}{2})^{t+1}) + \sum_{\tau=t+1}^{\infty} p^{(T+\tau)}x_{T+\tau}^{0}$ for all $n \ge n_t$.

By hypothesis, $\sum_{t=0}^{\infty} |p^{(t)}| < +\infty$ and for all (t, n)∈ $\mathbb{N}^0\times\mathbb{N}$, $x_{T+t}^{(n)} \in [0, b]$.

Thus, for all t∈ℕ and $n \ge n_t$, $\sum_{\tau=t+1}^{\infty} |p^{(T+\tau)}x_{T+\tau}^{N_t(n)}| = \sum_{\tau=t+1}^{\infty} |p^{(T+\tau)}|x_{T+\tau}^{N_t(n)}$.

For all τ∈ℕ and n∈ℕ, $x_\tau^{(n)} \in [0, b]$ <u>implies</u> for all τ∈ℕ and τ ≥ t+1, $\sup_{n \ge n_t} x_{T+\tau}^{N_\tau(n)}$ exists, $\sup_{n \ge n_t} x_{T+\tau}^{N_\tau(n)} \in [0, b]$ and $\sup_{n \ge n_t} x_{T+\tau}^{N_\tau(n)} \ge x_{T+t}^{N_\tau(n)}$ for all t ≥ τ+1 and $n \ge n_\tau$.

Thus, for all t∈ℕ and $n \ge n_t$, $\sum_{\tau=t+1}^{\infty}|p^{(T+\tau)}|[\sup_{n \ge n_t} x_{T+\tau}^{N_t(n)}] \ge \sum_{\tau=t+1}^{\infty} |p^{(T+\tau)}|x_{T+\tau}^{N_t(n)}$.

Further, for all t∈ℕ and $n \ge n_t$, $\sum_{\tau=t+1}^{\infty} |p^{(T+\tau)}x_{T+\tau}^{N_t(n)}| \ge \sum_{\tau=t+1}^{\infty} p^{(T+\tau)}x_{T+\tau}^{N_t(n)}$.

Thus, for all t∈ℕ, $V^T(x) + \sum_{\tau=t+1}^{\infty}|p^{(T+\tau)}|[\sup_{n \ge n_t} x_{T+\tau}^{N_t(n)}] \ge V^T(x) + \sum_{\tau=t+1}^{\infty} |p^{(T+\tau)}|x_{T+\tau}^{N_t(n)} \ge V^T(x) + \sum_{\tau=t+1}^{\infty} |p^{(T+\tau)}x_{T+\tau}^{N_t(n)}| \ge V^T(x) + \sum_{\tau=\tau+1}^{\infty} p^{(T+\tau)}x_{T+\tau}^{N_t(n)} > V^T(x^{(N_t(n))}) - \frac{\varepsilon}{4}(1- (\frac{1}{2})^{t+1}) + \sum_{\tau=t+1}^{\infty} p^{(T+\tau)}x_{T+\tau}^{0}$ for all $n \ge n_t$.

Thus, for all t∈ℕ, $V^T(x) + \sum_{\tau=t+1}^{\infty}|p^{(T+\tau)}|[\sup_{n \ge n_t} x_{T+\tau}^{N_t(n)}] > V^T(x^{(N_\tau(n))}) - \frac{\varepsilon}{4}(1- (\frac{1}{2})^{t+1}) + \sum_{\tau=t+1}^{\infty} p^{(T+\tau)}x_{T+\tau}^{0}$ for all $n \ge n_\tau$.

Thus, for all t∈ℕ, $V^T(x) + \sum_{\tau=t+1}^{\infty}|p^{(T+\tau)}|[\sup_{n \ge n_t} x_{T+\tau}^{N_t(n)}] \ge \limsup_{n\to\infty} V^T(x^{N_t(n)}) - \frac{\varepsilon}{4}(1- (\frac{1}{2})^{t+1}) + \sum_{\tau=t+1}^{\infty} p^{(T+\tau)}x_{T+\tau}^{0}$

By hypothesis, $\sum_{t=0}^{\infty} |p^{(t)}| < +\infty$ and for all (t, n)∈$\mathbb{N}^0\times\mathbb{N}$, $x_{T+t}^{(n)} \in [0, b]$.

Thus, $\lim_{t\to\infty} \sum_{\tau=t+1}^{\infty}|p^{(T+\tau)}|[\sup_{n \ge n_t} x_{T+\tau}^{N_t(n)}] = 0$.

Also, $\lim_{t\to\infty} \frac{\varepsilon}{4}(1 - (\frac{1}{2})^{t+1}) = \frac{\varepsilon}{4}$ and $\lim_{t\to\infty} \sum_{\tau=t+1}^{\infty} p^{(T+\tau)}x_{T+\tau}^{0} = 0$.

Thus, $V^T(x) = V^T(x) + \lim_{t\to\infty} \sum_{\tau=t+1}^{\infty}|p^{(T+t)}|[\sup_{n \ge n_t} x_{T+\tau}^{N_t(n)}] \ge \lim_{t\to\infty} \limsup_{n\to\infty} V^T(x^{N_t(n)}) - \lim_{t\to\infty} \frac{\varepsilon}{8}(1 - (\frac{1}{2})^{t+1}) + \lim_{t\to\infty} \sum_{\tau=t+1}^{\infty} p^{(T+\tau)}x_{T+\tau}^{0} = \lim_{t\to\infty} \limsup_{n\to\infty} V^T(x^{N_t(n)}) - \frac{\varepsilon}{4}$.

Thus, $V^T(x) \geq \lim_{t\to\infty} \limsup_{n\to\infty} V^T(x^{N_\tau(n)}) - \frac{\varepsilon}{4}$.

The above being true for all $\varepsilon > 0$, we get $V^T(x) \geq \lim_{t\to\infty} \limsup_{n\to\infty} V^T(x^{N_t(n)})$.

However, $\lim_{n\to\infty} V^T(x^{(n)})$ exists implies $\lim_{n\to\infty} V^T(x^{(n)}) = \limsup_{n\to\infty} V^T(x^{N_t(n)})$ for all $t\in\mathbb{N}$.

Thus, $V^T(x) \geq \lim_{t\to\infty} \limsup_{n\to\infty} V^T(x^{N_t(n)}) = V^T(x) \geq \lim_{t\to\infty} \lim_{n\to\infty} V^T(x^{(n)}) = \lim_{n\to\infty} V^T(x^{(n)})$, i.e., $V^T(x) \geq \lim_{n\to\infty} V^T(x^{(n)})$.

This combined with the continuity of $V^T$ on (0, b), $\lim_{x\to b} V^T(x) \geq V^T(b)$ and $\lim_{x\to 0} V^T(x) \geq V^T(0)$ implies $\lim_{x\to b} V^T(x) = V^T(b)$, $\lim_{x\to 0} V^T(x) = V^T(0)$ and hence the continuity of $V^T$ on X. Q.E.D.

The following result is analogous to the main results on dynamic programming available in Mitra (2000) and Sorger (2015)

**Proposition 5.2:** (i) For all $T\in\mathbb{N}^0$, $<V^t|\ t \geq T>$ satisfies the following ***functional equation of dynamic programming***: For all $x\in X$, $<x_t|\ t \geq T>\in\mathcal{S}^T(x)$ and $t \geq T$: $V^t(x_t) = p^{(t)}x_t + V^{t+1}(x_{t+1}) = p^{(t)}x_t + \max_{y\in\Omega_t(x_t)}\{V^{t+1}(y)\}$.

(ii) For all $x\in X$: $<x_t|t\in\mathbb{N}^0>\in\mathcal{S}(x)$ <u>if and only if</u> $<x_t|t\in\mathbb{N}^0>\in\mathcal{F}(x)$ and for all $T\in\mathbb{N}^0$ it is the case that $V^T(x_T) = p^{(T)}x_T + V^{T+1}(x_{T+1})$.

**Proof:** (i) Let $T\in\mathbb{N}^0$ and $x\in X$. Let $<x_t|\ t \geq T>\in\mathcal{S}^T(x)$.

Thus, $V^T(x) = V^T(x_T) = \sum_{\tau=T}^{\infty} p^{(\tau)}x_\tau$

Towards a contradiction suppose there exists $t \geq T$, such that $V^t(x_t) > p^{(t)}x_t + V^{t+1}(x_{t+1})$.

Thus, there exists $<y_{t+\tau}|\ \tau\in\mathbb{N}^0>\in\mathcal{F}^t(x_t)$ such that $y_t = x_t$ and $p^{(t)}x_t + \sum_{\tau=t+1}^{\infty} p^{(\tau)}y_\tau > p^{(t)}x_t + V^{t+1}(x_{t+1}) \geq p^{(t)}x_t + \sum_{\tau=t+1}^{\infty} p^{(\tau)}x_\tau$.

Hence, $\sum_{\tau=t+1}^{\infty} p^{(\tau)}y_\tau > V^{t+1}(x_{t+1}) \geq \sum_{\tau=t+1}^{\infty} p^{(\tau)}x_\tau$.

Clearly, $\{\tau \geq t|\ y_\tau \neq x_\tau\} = \{\tau \geq t+1|\ y_\tau \neq x_\tau\} \neq \phi$, since $<y_{t+\tau}|\ \tau\in\mathbb{N}^0>\in\mathcal{F}^t(x_t)$.

Let $<z_\tau|\ \tau \geq T>$ be such that $z_\tau = x_\tau$ for $\tau = T, \ldots, t$, $z_\tau = y_\tau$ for all $\tau \geq t+1$.

Thus, $<z_\tau|\ \tau \geq T>\in\mathcal{F}^T(x)$ and $\sum_{\tau=T}^{\infty} p^{(\tau)}z_\tau = \sum_{\tau=T}^{t} p^{(\tau)}z_\tau + \sum_{\tau=t+1}^{\infty} p^{(\tau)}z_\tau = \sum_{\tau=T}^{t} p^{(\tau)}x_\tau + \sum_{\tau=t+1}^{\infty} p^{(\tau)}y_\tau$.

Since, $\sum_{\tau=t+1}^{\infty} p^{(\tau)}y_\tau > \sum_{\tau=t+1}^{\infty} p^{(\tau)}x_\tau$, it follows that $\sum_{\tau=T}^{\infty} p^{(\tau)}z_\tau = \sum_{\tau=T}^{t} p^{(\tau)}x_\tau + \sum_{\tau=t+1}^{\infty} p^{(\tau)}y_\tau > \sum_{\tau=T}^{t} p^{(\tau)}x_\tau + \sum_{\tau=t+1}^{\infty} p^{(\tau)}x_\tau = \sum_{\tau=T}^{\infty} p^{(\tau)}x_\tau = V^T(x)$, i.e., $\sum_{\tau=T}^{\infty} p^{(\tau)}z_\tau > V^T(x)$, which is not possible, since $<z_\tau|\ \tau \geq T>\in\mathcal{F}^T(x)$.

Thus, it must be the case that $V^t(x) = p^{(t)}x_t + V^{(t+1)}(x_{t+1}) = p^{(T)}x_t + \max_{y\in\Omega_t(x_t)}\{V^{T+1}(y)\}$.

Thus, $<V^t|\ t \geq T>$ satisfies the "functional equation of dynamic programming".

(ii) Suppose $<x_t|\mathbb{N}^0>\in\mathcal{S}(x)$. Then from (i) it follows that $<x_t|t\in\mathbb{N}^0>\in\mathcal{F}(x)$ and for all $T\in\mathbb{N}^0$ it is the case that $V^T(x_T) = p^{(T)}x_T + V^{T+1}(x_{T+1})$.

Now suppose $<x_t|t\in\mathbb{N}^0>\in\mathcal{F}(x)$ and for all $T\in\mathbb{N}^0$: $V^T(x_T) = p^{(T)}x_T + V^{T+1}(x_{T+1})$.

Hence, $V(x) = V^0(x) = p^{(0)}x + V^1(x_1) = p^{(0)}x + p^{(1)}x_1 + V^2(x_2) = \ldots = \sum_{t=0}^{T} p^{(t)}x_t + V^{T+1}(x_{T+1})$ for all $T\in\mathbb{N}$.

Since, $|V^T(x_T)| \leq b\sum_{t=T}^{\infty}|p^{(t)}|$ and since $\lim_{T\to\infty}\sum_{t=T}^{\infty}|p^{(t)}| = 0$, it follows that $\lim_{T\to\infty} V^{T+1}(x_{T+1}) = 0$.

Thus, $V(x) = \lim_{T\to\infty}\sum_{t=0}^{T} p^{(t)}x_t + \lim_{T\to\infty} V^{T+1}(x_{T+1}) = \sum_{t=0}^{\infty} p^{(t)}x_t$.

Thus, $<x_t|t\in\mathbb{N}^0>\in \mathcal{S}(x)$. Q.E.D.

**6. Optimal Decision Rules:**

For all $T\in\mathbb{N}^0$ the correspondence $h^T: X\to\to X$ defined by $h^T(x) = \underset{y\in\Omega_T(x)}{\operatorname{argmax}} V^{(T+1)}(y)$ is said to be the **optimal period-T decision rule**.

As in Debreu (1959), $h^T$ will be said to be "upper semi-continuous" at $x\in X$ if for all $y\in h^T(x)$, sequences $<x^{(n)}|n\in\mathbb{N}>$, $<y^{(n)}|n\in\mathbb{N}>$ converging to x and y respectively: $[y^{(n)}\in h^T(x^{(n)})$ for all $n\in\mathbb{N}]$ implies $[y\in h^T(x)]$.

**Proposition 6.1:** Suppose that for all $T\in\mathbb{N}^0$, $\Omega_T$ is convex. Then:

(i) For all $T\in\mathbb{N}^0$, the optimal T-period rule is upper semi-continuous.

(ii) For all $T\in\mathbb{N}^0$ and $(x, y)\in\Omega_T$ and: $[y\notin h^T(x)]$ implies $[V^T(x) > p^{(T)}x + V^{T+1}(y)]$.

(iii) $<x_t|t\in\mathbb{N}^0>\in\mathcal{S}(x)$ <u>if and only if</u> $[x_0 = x$ and $x_{T+1}\in h^T(x_T)$ for all $T\in\mathbb{N}^0]$.

**Proof:** (i) Follows immediately from the definition of $h^T$, the assumption that $\Omega_T$ is a non-empty, closed and convex set and the continuity of $V^{T+1}$ on the closed and bounded interval $\{x\in|\Omega_T(x) \neq \phi\}$ established in part (ii) of proposition 4.1.

(ii) Follows immediately from the definition of $h^T$ and the functional equation of dynamic programming.

(iii) Follows immediately from the definition of h and part (ii) of proposition 5.2. Q.E.D.

**7. Linear Cake-eating Problems:**

Recall that the LDO problem is a linear cake-eating problem if $\Omega_t = \Omega = \{(z, y)|\ y\in[0, z], z\in[0, b]\}$ for all $t\in\mathbb{N}^0$.

**Note 12:** An interesting characteristic of a cake-eating problem is that for all $(x, T)\in X\times\mathbb{N}^0$: $[<x_t|\ t \geq T>\in \mathcal{F}^T(x)]$ <u>if and only if</u> $[x_T = x$ and the sequence $<x_t|\ t \geq T>$ is "non-increasing"].

In the cake-eating problem, given $x\in X$, the problem faced by the decision-maker is to choose $<x_t|t\in\mathbb{N}^0>$ so as to maximize $\sum_{t=\infty}^{\infty}p^{(t)}x_t$, subject to $(x_t, x_{t+1})\in\Omega$ for all $t\in\mathbb{N}^0$ and $x_0 = x$.

Note that if $p^{(t)} > 0$ for all $t\in\mathbb{N}^0$, then given $x\in X$, $<x_t|t\in\mathbb{N}^0>\in\mathcal{S}(x)$ for the cake-eating problem <u>if and only if</u> $x_t = x$ for all $t\in\mathbb{N}^0$, thereby implying the decision-maker never consumes any cake. Since the consumption of cake yields satisfaction, such a trajectory starting at x, would yield no satisfaction to the decision-maker.

On the other hand, if $p^{(t)} < 0$ for all $t \in \mathbb{N}^0$, then given $x \in X$, $\langle x_t | t \in \mathbb{N}^0 \rangle \in \mathcal{S}(x)$ for the cake-eating problem <u>if and only if</u> $x_{t+1} = 0$ for all $t \in \mathbb{N}^0$, thereby implying the decision-maker consumes the entire cake in period 0.

Thus, we need to be cautious in our interpretation of $p^{(t)}$ for $t \in \mathbb{N}$, although $p^{(0)} > 0$ is quite realistic and an innocuous assumption in this context.

Suppose the LDO problem is a linear cake-eating problem and for $T \in \mathbb{N}^0$ consider the following problem.

Given $x \in X$: Maximize $\sum_{t=T}^{\infty} p^{(t)} x_t$, subject to the infinite sequence $\langle x_t | \, t \geq T \rangle$ satisfying the constraints $(x_t, x_{t+1}) \in \Omega$ for all $t \geq T$ and $x_T = x$.

Clearly, $V^T(0) = 0$ for all $T \in \mathbb{N}^0$.

By proposition 5.1, we may conclude that for all $T \in \mathbb{N}^0$, $V^T$ is concave and continuous on X.

We have already noted that for $(x, T) \in X \times \mathbb{N}^0$: $\langle x_t | \, t \geq T \rangle \in \mathcal{F}^T(x)$ <u>if and only if</u> $x_T = x$ and $x_t \geq x_{t+1}$ for all $t \geq T$.

**Proposition 7.1:** Suppose that the LDO problem is a linear cake-eating problem. Let $T \in \mathbb{N}^0$.

(i) If $p^{(T)} \geq 0$, then for all x, y $\in X$, $y \geq x$ implies $V^T(y) \geq V^T(x)$ and thus $V^T(b) \geq V^T(x) \geq V^T(0)$ for all $x \in X$.

(ii) If $p^{(T)} > 0$, then for all x, $y \in X$, $y > x$ implies $V^T(y) > V^T(x)$.

(iii) Suppose $p^{(T)} \leq 0$. If for x, $y \in X$ such that $y \geq x$ there exists $\langle y_t | t \geq T \rangle \in \mathcal{S}^T(y)$ satisfying $y_{T+1} \leq x$, then $V^T(y) \leq V^T(x)$ and thus $V^T(b) \leq V^T(x) \geq V^T(0)$ for all $x \in X$.

(iv) Suppose $p^{(T)} < 0$. If for x, $y \in X$ such that $y > x$ there exists $\langle y_t | t \geq T \rangle \in \mathcal{S}^T(y)$ satisfying $y_{T+1} \leq x$, then $V^T(y) < V^T(x)$

(v) $D^-V^T(x)$ and $D^+V^T(x)$ exist for all $x \in (0, b)$, both $D^-V^T(x)$ and $D^+V^T(x)$ are non-increasing for $x \in (0, b)$ and for all $x \in (0, b)$, $V^T(y) - D^-V^T(x)(y) \leq V^T(x) - D^-V^T(x)x$ for all $y \in X$.

(vi) If $p^{(T)} \geq 0$ then $D^-V^T(b)$ exists $D^-V^T(b) \geq 0$, $D^-V^T(x) \geq D^-V^T(b)$ for all $x \in (a, b]$ and $V^T(y) - D^-V^T(b)(y) \leq V^T(b) - D^-V^T(b)b$ for all $y \in X$.

**Proof:** (i) Suppose $p^{(T)} \geq 0$. Let x, $y \in X = [0,b]$ with $y \geq x$ and let $\langle x_t | \, t \geq T \rangle \in \mathcal{S}^T(x)$. Let $\langle y_t | \, t \geq T \rangle$ be such that $y_T = y$ and $y_t = x_t$ for all $t \in \mathbb{N}$, $t > T$. Then, $\langle y_t | \, t \geq T \rangle \in \mathcal{F}^T(y)$.

Further, $\sum_{t=T}^{\infty} p^{(t)} y_t = p^{(T)}(y - x) + \sum_{t=T}^{\infty} p^{(t)} x_t \geq \sum_{t=T}^{\infty} p^{(t)} x_t$, since $p^{(T)} \geq 0$ and $y \geq x$ implies $p^{(T)} (y - x) \geq 0$.

Thus, $V^T(y) \geq V^T(x)$.

Thus, $V^T(b) \geq V^T(x) \geq V^T(0)$ for all $x \in X$.

(ii) If $p^{(T)} > 0$ and y > x, then $p^{(T)}$(y- x) > 0. Thus, if $<x_t|t \geq T> \in \mathcal{S}^T(x)$, then $\sum_{t=T}^{\infty} p^{(t)} y_t$ = $p^{(T)}$(y- x) + $\sum_{t=T}^{\infty} p^{(t)} x_t > \sum_{t=T}^{\infty} p^{(t)} x_t$ = $V^T(x)$.

Thus, $V^T(y) > V^T(x)$.

(iii) Let x, y∈X = [0,b] with y ≥ x and let $<y_t|t \geq T> \in \mathcal{S}^T(y)$.

Thus, $V^T(y) = \sum_{t=T}^{\infty} p^{(t)} y_t$.

Let $<z_t|t \geq T>$ be such that $z_T = x$ and $z_t = y_t$ for t > T.

$<z_t|t \geq T> \in \mathcal{F}^T(x)$ since $x \geq y_{T+1}$.

Thus, $\sum_{t=T}^{\infty} p^{(t)} z_t = \sum_{t=T}^{\infty} p^{(t)} y_t + p^{(T)}(x-y) \geq \sum_{t=T}^{\infty} p^{(t)} y_t$ since $p^{(T)} \leq 0$ and x ≤ y implies $p^{(T)}(x-y) \geq 0$.

Thus, $V(x) \geq \sum_{t=T}^{\infty} p^{(t)} z_t \geq \sum_{t=T}^{\infty} p^{(t)} y_t = V(y)$.

(iv) In (iii) $p^{(T)} < 0$ and x < y implies $p^{(T)}(x-y) > 0$, so that $V(x) \geq \sum_{t=T}^{\infty} p^{(t)} z_t > \sum_{t=T}^{\infty} p^{(t)} y_t$ = V(y).

(v) Follows from propositions 1 and 2 in Lahiri (2025) and from proposition 5.1 above. Proposition 5.1 says that $V^T$ is continuous and concave on X.

(vi) Follows from (i) of this proposition and part (ii) of proposition 3 in Lahiri (2025). Q.E.D.

**Note 12:** $p^{(T)} \geq 0$ is a sufficient condition for $V^{(T)}$ to be monotonically non-decreasing on X. $p^{(T)} \geq$ is not a necessary condition for $V^{(T)}$ to be monotonically non-decreasing on X. Hence it is quite possible that $<p^{(t)}| t \in \mathbb{N}^0>$ is an alternating sequence of positive and negative real numbers with $p^{(0)} > 0$ and yet for all $T \in \mathbb{N}^0$, $V^{(T)}$ is monotonically non-decreasing on X. However, if $V^{(T)}$ is monotonically non-decreasing, then some consequences follow immediately from such an assumption.

**Proposition 7.2:** Suppose the LDO problem is a linear cake-eating problem and let $T \in \mathbb{N}^0$. For x∈X, let $<x_t^*(x)| t \geq T>$ be such that for all t ≥ T, $x_t^*(x) = x$.

(i) If $V^{T+1}$ is non-decreasing, then for all x∈X, $x \in h^T(x)$.

(ii) If $V^{T+1}$ is strictly increasing, then for all x∈X, $h^T(x) = \{x\}$.

(iii) If $V^t$ is non-decreasing for all t > T, then for all x∈X, $<x_t^*(x)| t \geq T> \in \mathcal{S}^T(x)$.

**Proof:** (i) If $V^{T+1}$ is non-decreasing then for x∈X, $h^T(x) = \underset{y \in [0,x]}{\operatorname{argmax}} V^{(T+1)}(y)$ implies $x \in h^T(x)$.

(ii) If $V^{T+1}$ is strictly increasing, then from (i) we get that for all x∈X we get $x \in h^T(x)$. Further since $V^{T+1}$ is strictly increasing, 0 ≤ y < x implies V(y) < V(x), so that $h^T(x) = \{x\}$.

(iii) Follows by applying (i) to all t ≥ T. Q. E. D.

Proposition 7.2 implies that in the context of the cake-eating problem, assuming all time-dependent optimal value functions to be monotonically non-decreasing may lead to uninteresting consequences.

We now establish a "conditional" monotonicity property for $h^T$- the optimal T-period decision rule.

**Proposition 7.3:** For a linear cake-eating problem, let $T\in\mathbb{N}^0$ and x, y∈X such that x > y. Suppose $<x_t| t \geq T>\in\mathcal{S}^T(x)$ and $<y_t| t \geq T>\in\mathcal{S}^T(y)$. If $y_{T+1}\in h^T(y)\backslash h^T(x)$ and $y \geq x_{T+1}$, then $x_{T+1} > y_{T+1}$.

**Proof:** Since $y_{T+1}\in h^T(y)\backslash h^T(x)$, $y_{T+1} \neq x_{T+1}$.

Towards a contradiction suppose $x_{T+1} < y_{T+1}$.

$V^T(x) = p^{(T)}x + V^{T+1}(x_{T+1})$ and $V^T(y) = p^{(T)}y + V^{T+1}(y_{T+1})$.

Since we are concerned with a cake-eating problem, x > y implies that if $<\xi_t| t \geq T>$ is such that $\xi_T = x$, $\xi_t = y_t$ for all $t \geq T + 1$, then $<\xi_t| t \geq T>\in \mathcal{F}^T(x)$.

$y_{T+1}\in h^T(y)\backslash h^T(x)$ implies, $p^{(T)}x + V^{T+1}(x_{T+1}) > p^{(T)}x + V^{T+1}(y_{T+1})$ and hence $V^{T+1}(x_{T+1}) > V^{T+1}(y_{T+1})$.

Since we are concerned with a cake-eating problem, $y \geq x_{T+1}$ implies that if $<z_t| t \geq T>$ is such that $z_T = y$, $z_t = x_t$ for all $t \geq T+ 1$, then $<z_t| t \geq T>\in \mathcal{F}^T(y)$.

Thus, $V^T(y) = p^{(T)}y + V^{(T+1)}(y_{T+1}) \geq p^{(T)}y + V^{T+1}(x_{T+1})$ and hence $V^{(T+1)}(y_{T+1}) \geq V^{T+1}(x_{T+1})$, contradicting $V^{T+1}(x_{T+1}) > V^{T+1}(y_{T+1})$ that we obtained earlier.

Hence, it must be the case that $x_{T+1} \geq y_{T+1}$. Q. E. D.

A linear cake-eating problem $<(p^{(t)}, \Omega_t)| t\in\mathbb{N}^0>$ is said to be a **two-phase linear cake-eating problem** if there exists $T^+\in\mathbb{N}^0$, $T^-\in\mathbb{N}$ with $T^+ < T^-$ such that $p^{(t)} > 0$ for all $t \leq T^+$ and $p^{(t)} < 0$ for all $t \geq T^-$.

For such problems we have the following proposition about optimal decision rules.

**Proposition 7.4:** For a two-phase linear cake-eating problem for all $T\in\mathbb{N}^0$ and x∈X: $h^T(x) = \{x\}$ if $T + 1 \leq T^+$ and $h^T(x) = \{0\}$ if $T + 1 \geq T^-$.

**Proof:** If x = 0, then clearly $h^T(0) = \{0\}$ for all $T\in\mathbb{N}^0$. Hence, suppose x > 0.

Let $T\in\mathbb{N}^0$ satisfy $T + 1 \leq T^+$ and for some y < x suppose that $y\in h^T(x)$. Thus, there exists a non-increasing sequence $<y_t| t \geq T>$ with $y_T = x$, $y_{T+1} = y$ such that $<y_t| t \geq T>\in\mathcal{S}^T(x)$.

Consider the sequence $<x_t| t \geq T>$ with $x_T = x$, $x_{T+1} = x$, $x_t = y_t$ for all $t > T+1$.

Clearly, $<x_t| t \geq T>\in \mathcal{F}^T(x)$ and $\sum_{t=T}^{\infty} p^{(t)}x_t$- $\sum_{t=T}^{\infty} p^{(t)}y_t = p^{(T+1)}(x - y) > 0$ contradicting $<y_t| t \geq T>\in\mathcal{S}^T(x)$.

Since $h^T(x) \subset [0, x]$, it must be that $h^T(x) = \{x\}$.

Now suppose, that for $T\in\mathbb{N}^0$ it is the case that $T + 1 \geq T^-$ and for some y > 0 suppose that $y\in h^T(x)$. Thus, there exists a non-increasing sequence $<y_t| t \geq T>$ with $y_T = x$, $y_{T+1} = y$ such that $<y_t| t \geq T>\in\mathcal{S}^T(x)$.

Note that $y_t \geq 0$ for all $t \geq T$.

Consider the sequence $<x_t| t \geq T>$ with $x_T = x$, $x_t = 0$ for all $t \geq T+1$.

Clearly, $<x_t| t \geq T> \in \mathcal{F}^T(x)$ and $\sum_{t=T}^{\infty} p^{(t)}x_t - \sum_{t=T}^{\infty} p^{(t)}y_t = p^{(T+1)}(0 - y) + \sum_{t=T+2}^{\infty} p^{(t)}(0 - y_t) = - p^{(T+1)}y - \sum_{t=T+2}^{\infty} p^{(t)}y_t > - \sum_{t=T+2}^{\infty} p^{(t)}y_t$, since $p^{(T+1)} < 0$ and $y > 0$ implies $- p^{(T+1)}y > 0$.

However, $- \sum_{t=T+2}^{\infty} p^{(t)}y_t \geq 0$ since $p^{(t)} < 0$ and $y_t \geq 0$ for all $t \geq T+2$.

Thus, $\sum_{t=T}^{\infty} p^{(t)}x_t - \sum_{t=T}^{\infty} p^{(t)}y_t > - \sum_{t=T+2}^{\infty} p^{(t)}y_t \geq 0$, i.e., $\sum_{t=T}^{\infty} p^{(t)}x_t > \sum_{t=T}^{\infty} p^{(t)}y_t$.

This contradicts $<y_t| t \geq T> \in \mathcal{S}^T(x)$.

Since $h^T(x) \subset [0, x]$, it must be that $h^T(x) = \{x\}$. Q.E.D.

It follows immediately from proposition 7.4 that in the case of a "two-phase linear cake-eating problem", and any $x \in X$, if $<x_t| t \geq T> \in \mathcal{S}^T(x)$ then it must be the case that $x_t = x$ for all $t \leq T^+$ and $x_t = 0$ for all $t \geq T^-$.

**8. Interlinked Linear Dynamic Optimization Problems:**

In this section we consider the ILDO problem (see note 8 in section 3) $<((p_t^{(t)}, p_t^{(t+1)}), \Omega_t)| t \in \mathbb{N}^0>$ defined thus:

Given $x \in X$, Maximize $\sum_{t=0}^{\infty}[p_t^{(t)}x_t + p_t^{(t+1)}x_{t+1}]$ subject to the infinite sequence $<x_t|t \in \mathbb{N}^0>$ satisfying the constraints: $(x_t, x_{t+1}) \in \Omega_t$, $t \in \mathbb{N}^0$, $x_0 = x$.

The underlying assumption for such problems is that there exists a sequence $<(p_t^{(t)}, p_t^{(t+1)})| t \in \mathbb{N}^0>$ of "pairs of real numbers" satisfying $\sum_{t=0}^{\infty} |p_t^{(t)}| < +\infty$ and $\sum_{t=0}^{\infty} |p_t^{(t+1)}| < +\infty$.

Recall that for such problems the sequence $<(p^{(t)}, \Omega_t)| t \in \mathbb{N}^0>$, where $p^{(0)} = p_0^{(0)}$ and for $t \in \mathbb{N}$, $p^{(t)} = p_t^{(t)} + p_{t-1}^{(t)}$ is said to be the canonical LDO problem associated with the ILDO problem $<((p_t^{(t)}, p_t^{(t+1)}), \Omega_t)| t \in \mathbb{N}^0>$.

Example 5.4 in Sorger (2015) shows why the assumptions $\sum_{t=0}^{\infty} |p_t^{(t)}| < +\infty$ and $\sum_{t=0}^{\infty} |p_t^{(t+1)}| < +\infty$ are important for the definition of ILDO problems. That example is concerned with an ILDO $<((p_t^{(t)}, p_t^{(t+1)}), \Omega_t)| t \in \mathbb{N}^0>$ where for all $t \in \mathbb{N}^0$, $p_t^{(t)} = 1$ and $p_t^{(t+1)} = -1$. Clearly, $\sum_{t=0}^{\infty} |p_t^{(t)}| = +\infty$ and $\sum_{t=0}^{\infty} |p_{t+1}^{(t)}| = +\infty$.

The canonical LDO associated with it satisfies $p^{(0)} = 1$ and $p^{(t)} = 0$ for all $t \in \mathbb{N}$.

Let $<x_t|t \in \mathbb{N}^0>$ be an infinite sequence in X that converges to $x^* \in X$, $\sum_{t=0}^{\infty}[p_t^{(t)}x_t + p_t^{(t+1)}x_{t+1}] = \lim_{T\to\infty} \sum_{t=0}^{T}[p_t^{(t)}x_t + p_t^{(t+1)}x_{t+1}] = \lim_{T\to\infty} \sum_{t=0}^{T}[x_t - x_{t+1}] = x_0 - \lim_{T\to\infty} x_{T+1} = x_0 - x^* \neq x_0$, unless $x^* = 0$.

However, $\sum_{t=0}^{\infty} p^{(t)}x_t = x_0$.

In general, for any ILDO $<((p_t^{(t)}, p_t^{(t+1)}), \Omega_t)| t \in \mathbb{N}^0>$ and any infinite sequence $<x_t|t \in \mathbb{N}^0>$ in X, $\sum_{t=0}^{T}[p_t^{(t)}x_t + p_t^{(t+1)}x_{t+1}] = p_0^{(0)}x_0 + \sum_{t=1}^{T}(p_t^{(t)} + p_{t-1}^{(t)})x_t + p_T^{(T+1)}x_{T+1}$ for all $T \in \mathbb{N}$.

However, $\sum_{t=0}^{\infty}[p_t^{(t)}x_t + p_t^{(t+1)}x_{t+1}] = \lim_{T\to\infty}\sum_{t=0}^{T}[p_t^{(t)}x_t + p_t^{(t+1)}x_{t+1}] = p_0^{(0)}x_0 + \lim_{T\to\infty}[\sum_{t=1}^{T}(p_t^{(t)} + p_{t-1}^{(t)})\, x_t + p_T^{(T+1)}x_{T+1}] = p_0^{(0)}x_0 + \lim_{T\to\infty}[\sum_{t=1}^{T}(p_t^{(t)} + p_{t-1}^{(t)})\, x_t$ only if $\lim_{T\to\infty} p_T^{(T+1)}x_{T+1} = 0$, unlike in example 5.4 of Sorger (2015). One way to guarantee that $\lim_{T\to\infty} p_T^{(T+1)}x_{T+1} = 0$ is to assume $\sum_{t=0}^{\infty}|p_t^{(t)}| < +\infty$ and $\sum_{t=0}^{\infty}|p_{t+1}^{(t)}| < +\infty$, which is what we have done here.

For all $x\in X$ and $<x_t|t\in\mathbb{N}^0>\in\mathcal{F}(x)$: $<x_t|t\in\mathbb{N}^0>\in\mathcal{S}(x)$ <u>if and only if</u> $<x_t|t\in\mathbb{N}^0>\in \underset{<y_t|t\in\mathbb{N}^0>\in\mathcal{F}(x)}{\text{argmax}} \sum_{t=0}^{\infty}[p_t^{(t)}x_t + p_t^{(t+1)}x_{t+1}]$.

For all $(x, T)\in X\times\mathbb{N}$ and $<x_{T+t}|t\in\mathbb{N}^0>\in\mathcal{F}^T(x)$: $p_{T-1}^{(T)}x + \sum_{t=T}^{\infty}[p_t^{(t)}x_t + p_t^{(t+1)}x_{t+1}] = \sum_{t=T}^{\infty} p^{(t)}x_t$.

Thus, for $(x, T)\in X\times\mathbb{N}$ and $<x_{T+t}|t\in\mathbb{N}^0>\in\mathcal{F}^T(x)$: $<x_{T+t}|t\in\mathbb{N}^0>\in\mathcal{S}^T(x)$ <u>if and only if</u> $<x_{T+t}|t\in\mathbb{N}^0>\in \underset{<y_{T+t}|t\in\mathbb{N}^0>\in\mathcal{F}^T(x)}{\text{argmax}} [p_{T-1}^{(T)}x + \sum_{t=T}^{\infty}(p_t^{(t)}y_t + p_t^{(t+1)}y_{t+1})] = \underset{<y_{T+t}|t\in\mathbb{N}^0>\in\mathcal{F}^T(x)}{\text{argmax}} \sum_{t=T}^{\infty}(p_t^{(t)}y_t + p_t^{(t+1)}y_{t+1})$.

For all $T\in\mathbb{N}^0$, let $W^T: X \to\mathbb{R}$ be the function such that for all $x\in X$, $W^T(x) = \sum_{t=T}^{\infty}(p_t^{(t)}x_t + p_t^{(t+1)}x_{t+1})$ for $<x_{T+t}|t\in\mathbb{N}^0>\in\mathcal{S}^T(x)$.

Thus, for all $x\in X$, $W^0(x) = V^0(x) = V(x)$ and for all $(x, T)\in X\times \mathbb{N}$, $W^T(x) = V^T(x) - p_{T-1}^{(T)}x$.

The following result is an immediate consequence of proposition 5.1 and definition of $W^T$ for $T\in\mathbb{N}^0$.

**Proposition 8.1:** Suppose $\Omega_T$ is convex for all $T\in\mathbb{N}^0$. Then, for all $T\in\mathbb{N}^0$, $W^T$ is concave and continuous on X.

The only significant difference between an ILDO and an LDO is a reformulation of proposition 5.2, and most importantly the functional equation of dynamic programming.

**Proposition 8.2:** (i) For all $T\in\mathbb{N}^0$, $<W^t|\ t \geq T>$ satisfies the following <u>***functional equation of dynamic programming***</u>: For all $x\in X$, $<x_{T+t}|\ t\in\mathbb{N}^0>\in\mathcal{S}^T(x)$ and $t \geq T$: $W^t(x_t) = p_t^{(t)}x_t + p_t^{(t+1)}x_{t+1} + W^{t+1}(x_{t+1}) = p_t^{(t)}x + \max_{y\in\Omega_T(x_t)}\{p_t^{(t+1)}y + W^{(t+1)}(y)\}$.

(ii) For all $x\in X$: $<x_t|t\in\mathbb{N}^0>\in\mathcal{S}(x)$ <u>if and only if</u> $<x_t|t\in\mathbb{N}^0>\in\mathcal{F}(x)$ and for all $T\in\mathbb{N}^0$ it is the case that $W^T(x_T) = p_T^{(T)}x_T + p_T^{(T+1)}x_{T+1} + W^{T+1}(x_{T+1})$.

**Proof:** (i) Let $T\in\mathbb{N}^0$, $x\in X$ and $<x_{T+t}|\ t\in\mathbb{N}^0>\in\mathcal{S}^T(x)$.

From part (i) of proposition 5.2 we know that for $t \geq T$: $V^t(x_t) = p^{(t)}x_t + V^{t+1}(x_{t+1}) = p^{(t)}x + \max_{y\in\Omega_T(x_t)}\{V^{t+1}(y)\}$.

Thus, for $T = 0$, $W^0(x) = V^0(x) = p_0^{(0)}x + V^1(x_1) = p_0^{(0)}x + p_1^{(1)}x_1 + W^1(x_1)$.

The second equality follows from (i) of proposition 5.2.

From (i) of proposition 5.2 we know that $V^{1}(x_1) = \max_{y\in\Omega_0(x)} V^{1}(y) = \max_{y\in\Omega_0(x)} \{p_0^{(1)}y + W^{1}(y)\}$.

Thus, $W^0(x) = p_0^{(0)}x + p_1^{(1)}x_1 + W^{1}(x_1) = p_0^{(0)}x + \max_{y\in\Omega_0(x)} \{p_0^{(1)}y + W^{1}(y)\}$.

For T > 0, $W^T(x) = V^T(x) - p_{T-1}^{(T)}x = p^{(T)}x + V^{T+1}(x_{T+1}) - p_{T-1}^{(T)}x = p^{(T)}x + W^{T+1}(x_{T+1}) + p_T^{(T+1)}x_{T+1} - p_{T-1}^{(T)}x = p_T^{(T)}x + p_T^{(T+1)}x_{T+1} + W^{T+1}(x_{T+1})$.

The second equality follows from (i) of proposition 5.2.

From (i) of proposition 5.2 we know that $V^T(x) = p^{(T)}x + \max_{y\in\Omega_T(x)} V^{T+1}(y) = p^{(T)}x + \max_{y\in\Omega_T(x)} \{p_T^{(T+1)}y + W^{T+1}(y)\}$.

Hence, $W^T(x) = V^T(x) - p_{T-1}^{(T)}x = p^{(T)}x + \max_{y\in\Omega_T(x)} \{p_T^{(T+1)}y + W^{T+1}(y)\} - p_{T-1}^{(T)}x = p_T^{(T)}x + \max_{y\in\Omega_T(x)} \{p_T^{(T+1)}y + W^{T+1}(y)\}$.

Thus, $W^T(x) = p_T^{(T)}x + p_T^{(T+1)}x_{T+1} + W^{T+1}(x_{T+1}) = p_T^{(T)}x + \max_{y\in\Omega_T(x)} \{p_T^{(T+1)}y + W^{T+1}(y)\}$.

Towards a contradiction suppose that for some $T\in\mathbb{N}^0$ and t > T, it is the case that $W^t(x_t) > p_t^{(t)}x_t + p_t^{(t+1)}x_{t+1} + W^{t+1}(x_{t+1})$.

From part (i) of proposition 5.2 we know that $V^t(x_t) = p^{(t)}x_t + V^{t+1}(x_{t+1}) = p^{(t)}x_t + \max_{y\in\Omega_t(x_t)} \{V^{t+1}(y)\}$.

Since $W^t(x_t) = V^t(x_t) - p_{t-1}^{(t)}x_t$ and $W^{t+1}(x_{t+1}) = V^{t+1}(x_{t+1}) - p_t^{(t+1)}x_{t+1}$, [$W^t(x_t) > p_t^{(t)}x_t + p_t^{(t+1)}x_{t+1} + W^{t+1}(x_{t+1})$] implies [$V^t(x_t) - p_{t-1}^{(t)}x_t > p_t^{(t)}x_t + p_t^{(t+1)}x_{t+1} + V^{t+1}(x_{t+1}) - p_t^{(t+1)}x_{t+1}$], i.e., $V^t(x_t) > (p_t^{(t)}+p_{t-1}^{(t)})x_t + V^{t+1}(x_{t+1}) = p^{(t)}x_t + V^{t+1}(x_{t+1})$, leading to a contradiction.

Thus, it must be the case that $W^t(x_t) = p_t^{(t)}x_t + p_t^{(t+1)}x_{t+1} + W^{t+1}(x_{t+1}) \geq W^t(x_t) > p_t^{(t)}x_t + p_t^{(t+1)}y + W^{t+1}(y)$ for all $y\in\Omega_t(x_t)$.

Thus, $\langle W^t | t \geq T\rangle$ satisfies the "functional equation of dynamic programming" as stated in the statement of this proposition.

(ii) From (ii) of proposition 5.2 we know that $\langle x_t|t\in\mathbb{N}^0\rangle\in\mathcal{S}(x)$ <u>if and only if</u> $\langle x_t|t\in\mathbb{N}^0\rangle\in\mathcal{F}(x)$ and for all $T\in\mathbb{N}^0$ it is the case that $V^T(x_T) = p^{(T)}x_T + V^{T+1}(x_{T+1})$.

Thus, $\langle x_t|t\in\mathbb{N}^0\rangle\in\mathcal{S}(x)$ <u>if and only if</u> $\langle x_t|t\in\mathbb{N}^0\rangle\in\mathcal{F}(x)$ and for all $T\in\mathbb{N}^0$, $p_{T-1}^{(T)}x_T + W^T(x_T) = p^{(T)}x_T + p_T^{(T+1)}x_{T+1} + W^{T+1}(x_{T+1})$.

Thus, $\langle x_t|t\in\mathbb{N}^0\rangle\in\mathcal{S}(x)$ <u>if and only if</u> $\langle x_t|t\in\mathbb{N}^0\rangle\in\mathcal{F}(x)$ and for all $T\in\mathbb{N}^0$, $W^T(x_T) = p_T^{(T)}x_T + p_T^{(T+1)}x_{T+1} + W^{T+1}(x_{T+1})$. Q.E.D.

**Acknowledgment:** I wish to thank Subhadip Chakrabarti for comments on the paper.